\documentclass[11pt]{article}
\RequirePackage{fix-cm}

\def\B{\color{black}}

\usepackage{amstext,amsmath,amsfonts,amssymb,amsbsy,bm,color,graphicx,subfigure} 
\usepackage{color}
\usepackage{subfigure} 
\usepackage{lineno}
\usepackage{ulem}
\usepackage{epsfig,epstopdf,amsfonts }

\usepackage[lined,boxed,commentsnumbered]{algorithm2e}
\RestyleAlgo{ruled}

\newif\ifpdf
\ifx\pdfoutput\undefined
   \pdffalse
\else
   \pdfoutput=1
   \pdftrue
\fi
\ifpdf
   \usepackage{graphicx}
   \usepackage{epstopdf}
\else
   \usepackage{graphicx}
\fi

\newcommand{\ul}{\underline}
\newtheorem{thm}{Theorem}[section]

\begin{document}

\newtheorem{teo}[section]{Theorem}
\newtheorem{cor}[section]{Corollary}
\newtheorem{pro}[section]{Proposition}

\newcommand{\distt}{\mathop{\mathrm{dist}}}
\newcommand{\supp}{\mathop{\mathrm{supp}}}
\newcommand{\neigh}{\mathop{\mathrm{neigh}}}
\newcommand{\Trunc}{\mathop{\mathrm{Trunc}}}
\newcommand{\trunc}{\mathop{\mathrm{trunc}}}
\def\bfl{\mbox{\boldmath$\lambda$}}

\title{ On the accuracy of interpolation based on single-layer artificial neural networks with a focus on defeating the Runge phenomenon\B 
}


\author{ Ferdinando~Auricchio\footnote{ Dipartimento di Ingegneria Civile e Architettura, Università degli Studi di Pavia,
 Italy; auricchio@unipv.it},
  Maria~Roberta~Belardo\footnote{Scuola Superiore Meridionale, Naples, Italy; mariaroberta.belardo-ssm@unina.it }, \\ 
 Francesco Calabr\`{o}\footnote{Dipartimento di Matematica e Applicazioni ``Renato Caccioppoli", Universit\`a degli Studi di Napoli ``Federico II", Italy; francesco.calabro@unina.it},  Gianluca Fabiani\footnote{Scuola Superiore Meridionale, Naples, Italy; gianluca.fabiani@unina.it } \B 
\&  Ariel F.~Pascaner\footnote{ Dipartimento di Ingegneria Civile e Architettura, Università degli Studi di Pavia,
 Italy;  ariel.pascaner@unipv.it} }


\maketitle

\begin{abstract}
{Artificial Neural Networks (ANNs) are a tool in approximation theory widely used to solve interpolation problems. 
In fact, ANNs can be  assimilated to functions since they take an input and return an output.
The structure of the specifically adopted network determines the underlying approximation space, while the form of the function is selected by fixing the parameters of the network.
In the present paper, we consider one-hidden layer ANNs with a feedforward architecture, also referred to as shallow or two-layer networks, so that the structure is determined by the number and types of neurons.
The determination of the parameters that define the function, called training, is done via the resolution of the approximation problem, so by imposing the interpolation through a set of specific nodes.
We present the case where the parameters are trained using a procedure that is referred to as Extreme Learning Machine (ELM) that leads to a linear interpolation problem.
In such hypotheses, the existence of an ANN interpolating function is guaranteed.
\\
 Given that the ANN is interpolating, the error incurred occurs outside the sampling interpolation nodes provided by the user. In this study, various choices of nodes are analyzed: equispaced, Chebychev, and randomly selected ones. Then, the focus is on regular target functions, for which it is known that interpolation can lead to spurious oscillations, a phenomenon that in the ANN literature is referred to as overfitting. \B
We obtain good accuracy of the ANN interpolating function in all tested cases using these different types of interpolating nodes and different types of neurons.
 The following study is conducted starting from the well-known bell-shaped Runge example, which makes it clear that the construction of a global interpolating polynomial is accurate only if trained on suitably chosen nodes, ad example the Chebychev ones. \B
In order to evaluate the behavior when the number of interpolation nodes increases, we increase the number of neurons in our network and compare it with the interpolating polynomial.
We test using Runge's function and other well-known examples with different regularities.
As expected, the accuracy of the approximation with a global polynomial increases only if the Chebychev nodes are considered.
Instead, the error for the ANN interpolating function always decays, and in most cases we observe that the convergence follows what is observed in the polynomial case on Chebychev nodes, despite the set of nodes used for training.  Then we can conclude that the use of such an ANN defeats the Runge phenomenon.\B
\\
Our results show the power of ANNs to achieve excellent approximations when interpolating regular functions also starting from uniform and random nodes, particularly for Runge's function.}\\
Keywords: Feedforward Neural Networks; Function approximation; Runge's function; Extreme Learning Machine\\
Subclass {65D05  }
\end{abstract}


\section{Introduction}
\label{sec:Introduction}

With no doubt, Artificial Neural Networks (ANNs) are a numerical analysis tool with a great versatility and, therefore, they have had a great development in recent years \cite{higham2019deep,Mis2021,weinan2020towards}.
ANNs are a class of machine learning algorithms that have a well-established structure, flexibility, and approximation properties.
For this reason, they are widely used for function approximations, classification, regression, feature selection, and further applications are constantly being explored.
Nevertheless, they are often considered as a black-box for its complex and nonlinear structure.
Moreover, ANNs have open issues on some aspects of efficiency and accuracy when the number of training parameters increases.
Furthermore, ANNs have opened the road to Artificial Intelligence with the possibility to develop Digital Twins when a warehouse of input-output data is available, even in cases where the governing laws are unknown, even up to being able to mimic extremely complex phenomena, such as brain reaction.

Of all the highlighted points, the present contribution focuses on the aspect that ANNs can be profitably used in function approximation due to the Universal Approximation (UA) results, even under very mild hypotheses see, e.g., \cite{kratsios2021universal}.
In particular, UA states that ANNs can approximate a target function with any required precision, i.e., it can be defined an ANN that is close enough to any target function (at least asymptotically) also in spaces endorsed with very low regularity.

Moreover, some ANNs avoid the phenomenon called \emph{curse of dimensionality}, an expression coined by \cite{bellman1957dimensionality}, referring to a set of issues that appears when dealing with high-dimensional data, that are not present when working on low-dimensional problems.
For instance, the amount of data needed to obtain a sufficiently good representation of a given domain often increases exponentially with the dimension of the data.
The fact that ANNs avoid the curse of dimensionality opens the possibility to achieve efficient approximations, even in high dimensional cases see, e.g., \cite{han2018solving,jagtap2021deep,jin2017deep,karniadakis2021physics,lu2021learning}.

The construction of ANNs involves the determination of a sufficiently large number of functions (referred to in this context as \emph{activation functions}) connected in a grid structure composed by nodes.
The activation functions and the interactions among nodes are defined with parameters, which are determined via conditions that are intended as an interpolation (or approximation) of the available information.
This process is called the \emph{training} of the network and this is done by \emph{tuning} the parameters.
The first step to train an ANN is to initialize all the parameters in a proper way and then the tuning is performed, usually via optimization strategies, such as the minimization of a so-called \emph{loss function} see, e.g., \cite{bishop2006pattern,cyr2020robust}.

The existence of a ``good" interpolation network is guaranteed by the UA results and the large number of available tuning parameters.
Such a network can then ``remember" the information with which it has been trained and, even more important, is able to give an appropriate answer to new cases, which had not been used on the training process.
Translated into mathematical approximation theory language, a trained network is capable to do usual interpolation and extrapolation procedures.
On the other hand, from the authors' point of view, one of the main drawbacks of such a procedure is related to the unclear mathematical interpretation of the underlying approximating space and the definition of the training procedure in presence of a (quasi-) optimal set of parameters  \cite{vidal2017mathematics,weinan2020towards}.

A well known issue in approximation theory is the poor behavior of the degree-increasing polynomial interpolation when applied on some specific functions.
In particular, a paradigmatic example of this poor behavior is Runge's function see, e.g., \cite{corless2020runge}; when taking uniformly or randomly distributed sampling points of Runge's function and performing a polynomial interpolation, the difference between approximating and target function outside from the sampling points can diverge as the polynomial order increases.
One way to avoid this divergence is to sample the function using {ad-hoc} nodes, such as Chebychev nodes.
However, in many real-world, non theoretical problems {ad-hoc} optimal sampling cannot be performed, either because data are already given, or because the analytical expression of the underlying function is unknown, or a combination of both.

According to the previous discussion, the main aim of the present paper is to give some insight into the ability of  a class of ANNs to approximate functions, and in particular, we focus on the Runge phenomenon, that is, the presence of spurious oscillations in the interpolant function. The present is not the first study on the interpolation properties of such ANNs and the authors have focused on different aspects (among others, we refer to \cite{qu2016two,wang2011study,yuan2011optimization}), but none, to our knowledge, has results on this phenomenon. 
Then, we report \B results for the Runge's function and, by following examples in \cite{trefethen2019approximation}, also for other benchmark functions that exhibit different regularities.
Our reference method is the polynomial approximation obtained by interpolation in Chebychev nodes, as extensively discussed in \cite{trefethen2019approximation}; such a method is well suited to achieve high degree approximation, as well as for its relationships to the truncation (or projection) of the Chebychev series, and it is widely used also in other applications, such as in numerical quadrature see, e.g. \cite{trefethen2008gauss,calabro2009evaluation}.
Moreover, a stable reference implementation is available via CHEBFUN routines in MATLAB \cite{battles2004extension,driscoll2014chebfun}, so that the comparison can be easily conducted.

The paper is organized as follows.
Section \ref{sec:Preliminaries} presents an introduction to ANNs focusing on the definitions of the main concepts and  reviewing available results, with a focus on the properties related to the interpolation.
In particular, we set the specific ANN structure, later adopted in the paper, namely a single hidden layer network, trained via Extreme Learning Machine with over-parametrization and solved by Least Squares.
In Section \ref{sec:NumericalResults} we show numerical results. First, we use Runge's function, where we also report
results in the square (not overparametrized) case and for which we evaluate both the approximation of the function and of the derivative, being the derivative of the ANN interpolating function trivial to compute.
The section continues with the results  relative to other examples taken from Trefethen (\cite{trefethen2019approximation}), using overparametrized networks and focusing on the approximation of the interpolated functions.
Lastly, in Section \ref{sec:Conclusion} we present the final thoughts, conclusions and possible future directions.

\section{Preliminaries on Neural Networks}
\label{sec:Preliminaries}
In this section we define the main concepts related to Artificial Neural Networks (ANNs).
An ANN is a set of connected nodes, usually called \emph{neurons}, where each neuron takes different real-valued inputs and produces a real-valued output.
ANNs may have different connection structures between the neurons, called the \emph{architecture of the ANN}; for instance, common architectures are: feedforward, deep feedforward, perceptron, convolutional, recurrent, autoencoder, fully complex.
In the present work, we consider the most common architecture: the \emph{feedforward network} (FFN), described in Section \ref{subsec:ANNs}.

The process that each neuron performs is mainly divided into two steps, with each single step being represented by a function.
The first step is called \emph{interaction scheme} and combines the multiple inputs received from the neurons to produce a scalar output.
In Section \ref{subsec:InteractionSchemes} we describe some of the possible choices for the interaction scheme.
Then, the second step is represented by a transfer function, usually called \emph{activation function}, which takes the output of the interaction scheme as an input and generates another real value, that is the output of the neuron itself: in Section \ref{subsec:ActivationFunctions} we describe some of the possible choices for the activation functions.
The interaction scheme determines how the information coming from a set of inputs is combined, whereas the activation function modulates the result.
As discussed in more detail in the respective sections, the choice of each of these steps determines the behavior of the ANN.

After these introductions of ANNs, in Section \ref{subsec:MathANN} we focus the analysis and particularise the equations to the case of a single hidden layer.
Finally, in Section \ref{subsec:ELM} we explain the adopted training procedure, called Extreme Learning Machine (ELM).

Throughout the paper we use the following notation.
Scalar values are represented by letters without underlining, whereas vectors and matrices are indicated with a single or a double underline, respectively.
The elements of a vector are represented with the same letter as the vector itself, with a single sub-script indicating the number of elements.
Similarly, the elements of a matrix are represented with the same letter as the matrix itself, using two sub-scripts.
In both cases, the letter with the sub-script(s) is not underlined, since it corresponds to a scalar value.
Moreover, a single row or column of a matrix (hence a vector) is denoted with the same letter as the matrix and a single sub-script.
In this case, the row or column is underlined with a single line since it corresponds to a vector.
Finally, super-scripts between parentheses always indicate the corresponding layer number within the ANN.


\subsection{Architecture of ANNs: the feedforward case} \label{subsec:ANNs}
The first choice for the determination of the structure of an ANN is the architecture of the grid.
We use the so-called feedforward architecture, resulting in the feedforward network (FFN), where neurons are stored in consecutive layers; the layers are indicated with an index $l$, with $l=0,1,\dots,\mathcal{L},\mathcal{L}+1$.
The FFN architecture is such that at each layer the input is taken from the output of the previous layer, elaborated by the interaction scheme and the activation functions of the neurons present in the considered layer, and passed as an input to the next layer.
The first layer of the FFN ($l=0$) is called \emph{ANN input layer} or simply \emph{input layer} and it is not actually composed by neurons; each unit of this first layer has the only role of reading the input data and to feed them to the following layer of neurons to start the network process; accordingly, the amount of units in the first layer is always equal to the number of inputs of the network.
The last layer of the FFN ($l=\mathcal{L}+1$) is called \emph{ANN output layer} or simply \emph{output layer} since its neurons' outputs constitute the output of the network itself.
All the layers between the ANN input and output layers ($l=1,\dots,\mathcal{L}$) are called \emph{hidden layers}.
Each $l$--th layer consist of $N^{(l)}$ neurons, each one with its own real scalar output.
The set of output values of the $l$--th layer's neurons\footnote{As mentioned, the units of the input layer do not perform any of the operations usually associated with the concept of a neuron. However, we refer to the units of all layers (including the input layer) as \emph{neurons} for consistency with the literature and for clarity and simplification of the mathematical expressions.}
can be thus thought as a vector, arranged, for example, in the form: $\ul{x}^{(l)}=(x_1^{(l)},\dots,x_{N^{(l)}}^{(l)}) \in \mathbb{R}^{N^{(l)}}$.
In the FFN architecture, each neuron of the $l$--th layer is fully connected with the nodes of the previous $(l-1)$--th layer by weighted edges (often called just \emph{weights}); the weights are scalar values associated to the connection between neurons of adjacent layers.

When the network has two or more hidden layers ($\mathcal{L} \ge 2$), the architecture is usually referred to as \emph{deep feedforward network (DFFN)}.
Instead, when the network has just one hidden layer ($\mathcal{L} = 1$), the architecture is usually referred to as \emph{non-deep FFN}, or \emph{shallow neural networks} \cite{pinkus,siegel2022high} or \emph{two-layer networks} \cite{ma2022barron,siegel2020approximation}.
This last terminology evidences the fact that the input layer does not perform any operation, hence such ANN relies only on a single hidden layer and the output layer.


\subsection{Interaction schemes}
\label{subsec:InteractionSchemes}

The first step performed by the neurons of the generic $l$--th layer is represented by the interaction scheme, which combines all the outputs of the neurons of the previous layer (i.e. the real values $x_j^{(l-1)}$, with $j=1\dots N^{(l-1)}$) to produce a single output value for each single neuron (i.e., the real values denoted in the following as $z_i^{(l)}$, with $i=1\dots N^{(l)}$).

For each $i$--th neuron on the $l$--th layer, such a combination of the output neurons of the previous layer is based on the use of a vector of weights, indicated as $\ul{A}_i^{(l)} \in \mathbb{R}^{N^{(l-1)}}$, and a scalar bias, indicated as $\beta_i^{(l)}\in\mathbb{R}$, with $i=1\dots N^{(l)}$. In general, the interaction scheme of the single neuron is represented as a function $\kappa_i^{(l)}:(\ul{x}^{(l-1)},\ul{A}_i^{(l)},\beta_i^{(l)}) \in \mathbb{R}^{N^{(l-1)}}\times \mathbb{R}^{N^{(l-1)}} \times \mathbb{R} \rightarrow z_i^{(l)} \in \mathbb{R}$, with $i=1\dots N^{(l)}$, so that:
\begin{equation}
     z_i^{(l)}=\kappa_i^{(l)}(\ul{x}^{(l-1)},\ul{A}_i^{(l)},\beta_i^{(l)}).
     \label{eq:interaction_scheme}
\end{equation}
Figure \ref{fig:neuron} depicts the interaction scheme of the generic $i$--th neuron on the $l$--th layer (the figure also shows the activation function $\psi_i^{(l)}$, which is addressed in Section \ref{subsec:ActivationFunctions}).

It is possible to organize the weights $\ul{A}_i^{(l)}$ and biases $\beta_i^{(l)}$ associated to the $i$--th neuron of the $l$--th layer into a matrix and a vector, respectively.
We define the weights matrix $\ul{\ul{A}}^{(l)}\in \mathbb{R}^{N^{(l)}\times N^{(l-1)}}$ by laying the vectors $\ul{A}_i^{(l)}$ as rows.
Likewise, the vector of biases $\ul{\beta}^{(l)}\in \mathbb{R}^{N^{(l)}}$ is defined organizing the individual biases as a vector.
Consequently, $\ul{\ul{A}}^{(l)}$ and $\ul{\beta}^{(l)}$ contain the weights and biases associated to all the neurons of the $l$--th layer.



The interaction scheme distinguishes different classes of FFN formulations.
In the following we address only the additive and distance-like schemes, referring to them as classic DFFN and RBF networks, respectively, but it is worthwhile to remark that ANNs have been also extended to many others schemes.

\subsubsection{Additive interaction scheme}
\label{subsubsec:AdditiveScheme}
Such a scheme is based on an additive combination of the weights $\ul{A}_i^{(l)}$, biases $\beta_i^{(l)}$ and inputs $\ul{x}^{(l-1)}$, and it is defined as follows:
\begin{equation}
     z_i^{(l)}=
     \sum_{j=1}^{N^{(l-1)}} A_{ij}^{(l)} x_j^{(l-1)}+\beta_i^{(l)}
     \qquad i=1,\dots,N^{(l)},
 \end{equation}
 where $A_{ij}^{(l)}$ is the $i$--th row and $j$--th column element of the matrix $\ul{\ul{A}}^{(l)}$.

\subsubsection{Distance-like interaction scheme}
\label{subsubsec:Distance-likeScheme}

Such a scheme is based on the extension of the notion of Radial Basis Functions (RBF) to the networks and it is defined as follows.

\begin{equation}
    z_i^{(l)}=
    \frac{||\ul{x}^{(l-1)}-\ul{A}_i^{(l)}||_2}{\beta_i^{(l)}}
    \qquad i=1,\dots,N^{(l)},
\end{equation}
where $\ul{A}_i^{(l)}$ is the $i$--th row of the matrix $\ul{\ul{A}}^{(l)}$.

{Since internal parameters $(\ul{A}_i^{(l)},\beta_i^{(l)})$ are difficult to train, one common choice is to fix the radius $\beta_i^{(l)}$ and set random centers $\ul{A}_i^{(l)}$ and then proceed to the solution of the linear system in the weights of the last layer $\ul{A}_i^{(\mathcal{L}+1)}$ see, e.g., \cite{broomhead1988radial,hryniowski2019deeplabnet}.
Such a procedure is similar to the one proposed in the training via ELMs, as we present in Section \ref{subsec:ELM}.}


\subsection{Activation functions}
\label{subsec:ActivationFunctions}
The second step performed by each neuron in the $l$--th layer is to modulate the output of the interaction scheme $z_i^{(l)}$ by means of the activation function to produce the so-called \emph{output} of the neuron, which is indicated as $x_i^{(l)}$.
The output $x_i^{(l)}$ is computed by the activation function $\psi_i^{(l)}:z^{(l)}_i \in \mathbb{R} \rightarrow x_i^{(l)} \in \mathbb{R}$ as:
\begin{equation}  \label{eq:interaction_scheme_2} 
    x_i^{(l)}=\psi_i^{(l)}(z_i^{(l)})
		\qquad i=1,\dots, N^{(l)}.
\end{equation}

\begin{figure}
\centering
\includegraphics[scale=0.4]{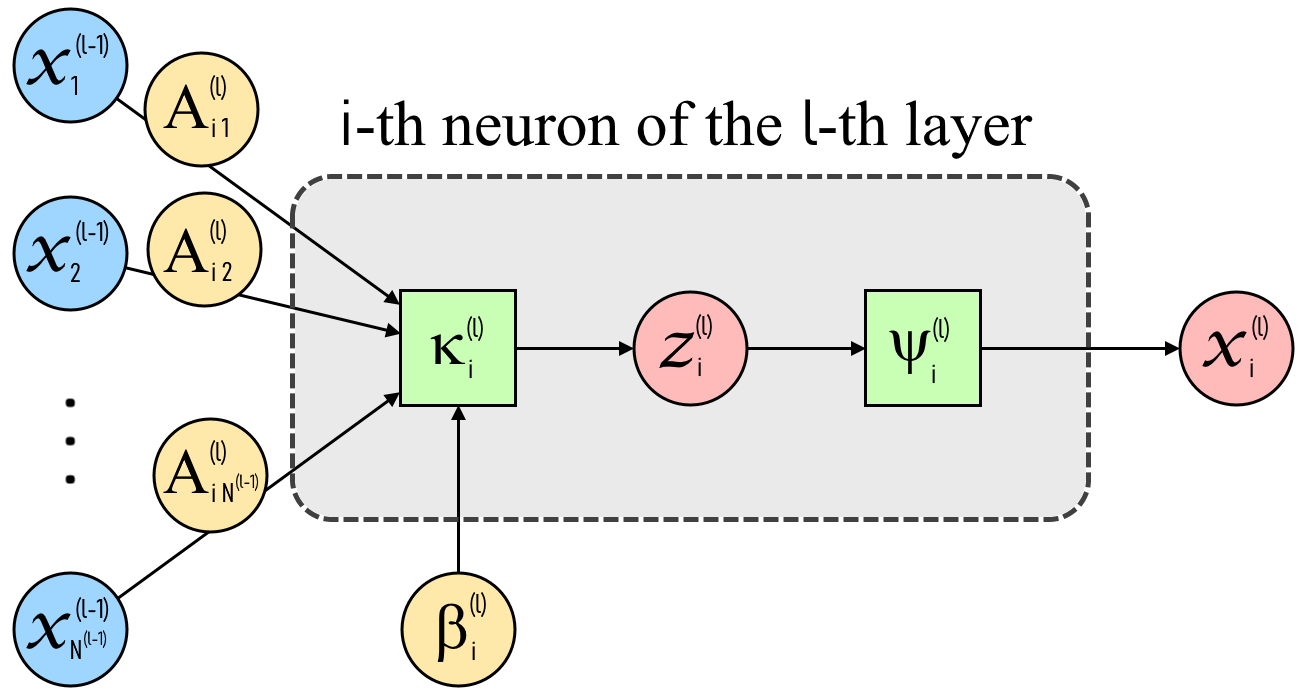}
\caption{Schematic representation of the action of the generic $i$--th neuron of the $l$--th layer. The inputs $x_j^{(l-1)}$ are represented with blue circles, the weights $A_{ij}^{(l)}$ and bias $\beta_i^{(l)}$ are represented by yellow circles, the computed values $z_i^{(l)}$ and $x_i^{(l)}$ are represented with red circles and the functions performed by the neuron (i.e. the interaction scheme $\kappa_i^{(l)}$ and the activation function $\psi_i^{(l)}$) are represented with green squares.}
\label{fig:neuron}
\end{figure}

Figure \ref{fig:neuron} shows a representation of the action of the generic $i$--th neuron on the $l$--th layer, including both the interaction scheme $\kappa_i^{(l)}$ and the activation function $\psi_i^{(l)}$.

The use of activation functions gives a wide range of possibilities when defining a network.
In principle, one can mix the use of different activation functions in the different hidden layers.
However, we use just one activation function for all neurons of all layers, i.e.: $\psi_i^{(l)}=\psi$.
In the following, we report the choices that are more common and that we consider in our tests.

\subsubsection{Activation functions used with additive interaction scheme}
The most common choices for activation functions adopted in combination with additive interaction scheme are: 
\begin{itemize}
    \item Logistic Sigmoid (LS):
    \[\psi(x)=\frac{1}{1+\text{exp}(-x)}\]
    \item SoftPlus (SP):
    \[\psi(x) = \log(1+\text{exp}(x)) . \]

\end{itemize}
LS and SP functions are widely used when performing approximations in unbounded domains or for extrapolating.

Another common activation function adopted in combination with additive interaction scheme is the Rectified Linear Unit (ReLU), which is defined as $\psi(x) = max\{x,0\}$. Notice that the SP function is a smooth $C^{\infty}$ approximation of the ReLU function with the property that the derivative of the SP function is the LS function.
It is well known that shallow networks with ReLU activation functions have piece-wise linear functions as an output.
Since our aim is to obtain regular functions on the output of the ANN, in our tests we do not implement networks with ReLU activation functions.

Finally, we recall that in the case of shallow networks with LS activation functions the network gives the formation of the so-called \emph{Ridge} functions \cite{pinkus2015ridge}, or superimpositions of LS functions  \cite{barron1993universal}.


\subsubsection{Activation functions used with distance-like interaction scheme}

A common choice for activation function adopted in combination with distance-like interaction scheme is:
\begin{itemize}
    \item Gaussian Radial Basis (GRB):
    \[\psi(x)=\text{exp}(-x^2)\]
\end{itemize}

{
GRB functions are widely used in different domains, and the understanding of the properties of these is well established,  see e.g. \cite{fornberg2011stable}.
}



\subsection{Mathematical formulation of the entire ANN} \label{subsec:MathANN}

According to the notation presented above, for each ANN layer, we can define a \emph{layer map} $\ul{\Phi}^{(l)}:\ul{x}^{(l-1)}\in\mathbb{R}^{N^{(l-1)}} \rightarrow \ul{x}^{(l)}\in\mathbb{R}^{N^{(l)}}$ for $l=1,\dots,\mathcal{L}+1$ as:
\begin{equation}
   \ul{x}^{(l)}= \ul{\Phi}^{(l)}(\ul{x}^{(l-1)})=\biggl(\psi_1^{(l)}( z_1^{(l)}),\dots,\psi_{N^{(l)}}^{(l)}(z_{N^{(l)}}^{(l)})\biggr)
   \qquad
    \textnormal{with }
    z_i^{(l)}=\kappa_i^{(l)}(\ul{x}^{(l-1)},\ul{A}_i^{(l)},\beta_i^{(l)})
   \label{eq:layermap}
\end{equation}

Finally, the whole ANN can be interpreted as a map $\mathcal{F}:\ul{x}^{(0)}\in \mathbb{R}^{N^{(0)}}\rightarrow \ul{x}^{(\mathcal{L}+1)} \in \mathbb{R}^{N^{(\mathcal{L}+1)}}$, i.e. a map between the input domain and the output domain that can be expressed by:
\begin{equation}\label{eq:composizione}
    \ul{x}^{(\mathcal{L}+1)}=\mathcal{F}(\ul{x}^{(0)})= \ul{\Phi}^{(\mathcal{L}+1)}\circ\dots \circ \ul{\Phi}^{(1)}(\ul{x}^{(0)})
\end{equation}
where the symbol $\circ$ is the function composition, \emph{i.e.}, $f \circ g (x)= f(g(x))$.

When FFNs are used for regression or approximation purposes in the context of continuous, integrable or square-integrable functions, a common choice is to use linear activation functions for the neurons on the last layer.
The final output is taken as a linear combination of the evaluations done at the last hidden layer, so that transfer functions of the last layer $\psi_i^{(\mathcal{L}+1)}$ are the identity functions:
\begin{equation}\label{eq:ultimolayer}
    x_i^{(\mathcal{L}+1)}=\sum_{j=1}^{N^{(\mathcal{L})}}A^{(\mathcal{L}+1)}_{ij} \psi_j^{(\mathcal{L})}(z_j^{(\mathcal{L})}) \qquad i=1,\dots,N^{(\mathcal{L}+1)}\, .
\end{equation}

In this case, the structure of the FFN is exactly a linear combination of the activation functions of the last hidden layer's neurons.
If such functions are chosen linearly independent, they are a basis of the approximating space, so this space is easily described.
{
Overall, the effect of the weights of the last hidden layer is only to combine the activation functions, whereas altering the internal parameters that define the shape of the functions changes the approximating space.
}

\subsubsection{Approximation results}
{
The universal approximation result can be enunciated for FFNs with the same activation function for alla the FFN as follows.
}

\begin{thm}\label{th:1}
    Consider a DFFN with the same activation function $\psi$ in all the neurons of all the layers. 
		Moreover let $\psi \in C(\mathbb{R})$. Then the space $\mathcal{M}^{(\mathcal{L})}(\psi)$ generated by a linear combination of the outputs of the last hidden layer $\ul{x}^{(\mathcal{L})}$, defined by:
    \begin{equation*}
        \mathcal{M}^{(\mathcal{L})}(\psi)=\sum_{i=1}^{N^{(\mathcal{L})}} w_i\psi(z_i^{(\mathcal{L})})
        \qquad
        \textnormal{with }
        z_i^{(\mathcal{L})}=\kappa_i^{(\mathcal{L})}(\ul{x}^{(\mathcal{L}-1)},\ul{A}_i^{(\mathcal{L})},\beta_i^{(\mathcal{L})})
    \end{equation*}
     is dense in $C(\mathbb{R})$, in the topology of uniform convergence on compacta, if and only if $\psi$ is not a polynomial.
\end{thm} 

Notice that the choice of the interaction scheme is not discussed here since the two cases presented in Section \ref{subsec:InteractionSchemes} fulfill the requested hypotheses.

We remark that the norm associated to the uniform convergence on compacta is very general, since if $\mu$ is a non--negative finite Borel measure defined on some compact set $K$, then $C(K)$ is dense in $L^p(K,\mu)$, for any $1\le p \le +\infty$.

{
The possibility of using different activation functions in different layers comes as a consequence of the universal approximation result, that states that a property on universal approximation can be obtained with every non-polynomial function see, e.g., \cite{leshno1993multilayer,hornik1990universal,park1991universal}.
}
In fact, Theorem \ref{th:1} can be easily extended to networks with different activation functions in different layers; assuming to deal with two different activation functions $\psi_1$ and $\psi_2$ in different layers, then the generated space is the direct sum of two spaces dense in $C(\mathbb{R})$: $\mathcal{M}^{(\mathcal{L})}(\psi_1,\psi_2)=\mathcal{M}^{(\mathcal{L})}(\psi_1) \bigoplus \mathcal{M}^{(\mathcal{L})}(\psi_2)$.
Furthermore, from a density point of view, there is no need to consider a deep ANN, because, at least for the univariate case, a single hidden layer is enough to achieve the universal approximation propriety.
Moreover, the construction of deep neural networks gives spaces of composite functions, see \eqref{eq:composizione}, and such composition is unnecessary for the interpolation or approximation of univariate functions due to Theorem \ref{th:1}.  For this reason, in our analysis we focus on the most simple architecture and take in account only few possible choices of activation functions. With this, we are not intending to assert that this is the optimal choice for all possible cases, but rather because even in these frameworks, we find what aligns with the purpose of this study.\B 


Furthermore, we would like to point out that function composition can be interpreted as a construction of a new feature set for the output layer.
In fact, when a point $\ul{x}^{(0)}$ is given as an input, this is transformed in an element in the $N^{(1)}$--dimensional space by the action of the activation functions of first hidden layer; then, the activation functions of the second hidden layer - or, in the case of one layer, the activation function - has to deal with the new point $\ul{x}^{(1)}$ and, so on, this happens in all the intermediate layers.
Accordingly, one common point of view to look at the behavior of ANNs is to think that the effect of all hidden layers is to encode the input in $N^{(\mathcal{L})}$ new coordinates (i.e., the outputs of the last hidden layer), acting thus as a filter.
Then, these $N^{(\mathcal{L})}$ values are linearly combined to give the output, when linear activation functions for the neurons on the last layer are used.
In other words, the network takes the input data and converts it into numbers that, linearly combined with a weight vector, determine the behavior of the network.



\subsubsection{Single hidden layer FFN in the univariate scalar case}
\label{subsec:SingleHiddenLayerFFN}

Our aim now is to particularise the formulae from previous sections to the case of a single hidden layer ANN and univariate scalar functions, which are the focus of the present work.

In the case of a single hidden layer ($\mathcal{L}=1$), the FFN is composed only by three layers, numbered $l=0$ (input layer), $l=1$ (hidden layer) and $l=2$ (output layer); then, there are only two matrices of weights $\ul{\ul{A}}^{(1)}$ and $\ul{\ul{A}}^{(2)}$ and two vectors of biases $\ul{\beta}^{(1)}$ and $\ul{\beta}^{(2)}$.

\begin{figure}
\centering
\includegraphics[scale=0.4]{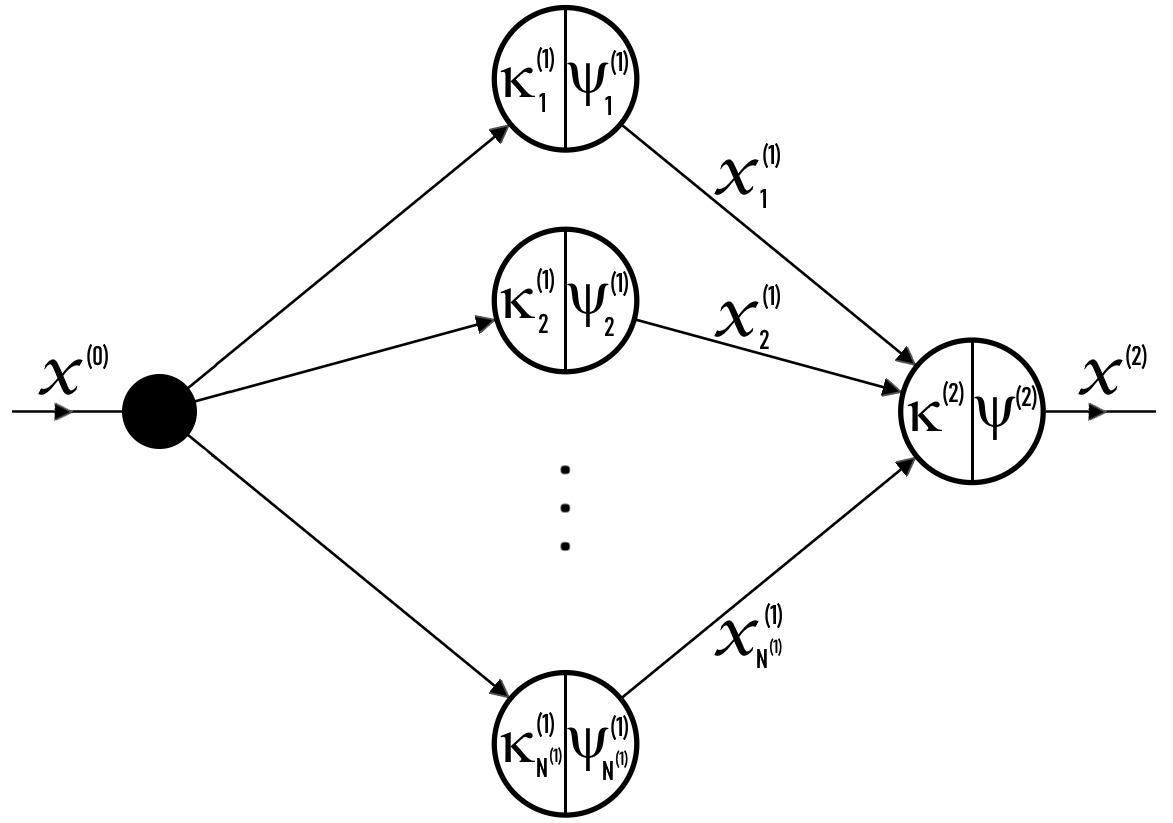}
\caption{Schematic representation of the single-hidden layer ANN with scalar input $x^{(0)}$ and scalar output $x^{(2)}$. Each neuron is represented by a circle divided into the interaction scheme $\kappa_i^{(l)}$ and the activation function $\psi_i^{(l)}$. The input neuron, which does not perform any of these operations, is represented by a black dot.}
\label{fig:network}
\end{figure}

Moreover, when the target functions to be interpolated are univariate and scalar (as we consider in our work), the ANN can be further simplified.
On the one hand, the fact that the functions are univariate implies that $N^{(0)}=1$, hence $x^{(0)}$ is a scalar and $\ul{A}^{(1)}$ reduces to a vector.
On the other hand, the fact that the functions are scalar implies that $N^{(2)}=1$, hence $x^{(2)}$ and $\beta^{(2)}$ reduce to scalars and $\ul{A}^{(2)}$ reduces to a vector, which we write as: $\ul{A}^{(2)}=\ul{w}=(w_1,w_2,\dots,w_{N^{(1)}}) \in \mathbb{R}^{N^{(1)}}$ and we denote as \emph{external weights}.
Figure \ref{fig:network} shows a diagram of the resulting ANN.

Under these conditions, equations (\ref{eq:composizione}) and (\ref{eq:ultimolayer}) become, respectively:

\begin{equation}\label{eq:composizione_SL}
\begin{split}
    & x^{(2)}=\mathcal{F}(x^{(0)})= \Phi^{(2)} (\ul{\Phi}^{(1)}(x^{(0)})),
    \\
    & \textnormal{where }
    \ul{\Phi}^{(1)}(x^{(0)})=\biggl(\psi_1^{(1)}( z_1^{(1)}),\dots,\psi_{N^{(1)}}^{(1)}(z_{N^{(1)}}^{(1)})\biggr)
    \textnormal{ and }
    \Phi^{(2)}\equiv\psi^{(2)}\,,\
\end{split}
\end{equation}

\begin{equation}\label{eq:ultimolayer_SL}
    x^{(2)}=\sum_{i=1}^{N^{(1)}}w_i \psi_i^{(1)}(z_i^{(1)})\ 
            \qquad
        \textnormal{with }
        z_i^{(1)}=\kappa_i^{(1)}(x^{(0)},A_i^{(1)},\beta_i^{(1)}).
\end{equation}


A classic FFN with a single hidden layer of $N^{(1)}$ neurons and with an additive interaction scheme can be written as:
\begin{equation}\label{eq:DFFN_1}
    x^{(2)}=\sum_{i=1}^{N^{(1)}}w_{i}\psi_i^{(1)}\biggl( A_{i}^{(1)} x^{(0)}+\beta_i^{(1)} \biggr) \ .
\end{equation}

Likewise, a RBF network with a single hidden layer of $N^{(1)}$ neurons can be written as:
\begin{equation}\label{eq:RBF_1}
    x^{(2)}=\sum_{i=1}^{N^{(1)}}w_{i}\psi_i^{(1)}\biggl(\frac{|{x}^{(0)}-{A}_i^{(1)}|}{\beta_i^{(1)}}\biggr) \ .
\end{equation}

\subsection{Training of the single-layer network}
\label{subsec:ELM}
When a choice has been done for the number of neurons $N^{(1)}$ and for the activation functions $\psi_i^{(1)}$ (see, for example, equations \eqref{eq:DFFN_1} and \eqref{eq:RBF_1}), then we are left with the determination of the unknown parameters $A_i^{(1)}$, $\beta_i^{(1)}$ and $w_i$, with $i=1\dots N^{(1)}$ (a process often indicated as training process).
Among all possible choices, we assume to fix a-priori all weights and biases associated to the neurons in the hidden layer, i.e. $A_i^{(1)}$ and $\beta_i^{(1)}$, choice often indicated under the name of \emph{Extreme Learning Machine} (ELM) and the ANNs trained this way are usually called \emph{Random Projection Networks}.
Accordingly, the only unknown of a FFN trained via ELM are the external weights $w_i$, i.e. the coefficients of the linear combination of the outputs of the hidden layer.  More detailed analysis on some parameters choices that can influence the analysis can be found, ad example, in \cite{dong2022computing,fornasier2022finite}. \B
The weights $w_i$ are selected throughout the optimization problem, that in our case is the satisfaction of the interpolation requirements.

Therefore, the number of degrees of freedom of an ELM network ($N_{dof}^{ELM}$) is equal to the number of neurons in the last hidden layer, independently of how many hidden layers are considered.
In our case:
\begin{equation}
    N_{dof}^{ELM}=N^{(1)},
\end{equation}
whereas for a generic classical DFFN, with $\mathcal{L}$ hidden layers, the number of degrees of freedom would be $N_{dof}^{DFFN}=N^{(\mathcal{L})}+\sum_{i=0}^{\mathcal{L}-1} (N^{(i)}+1)\cdot N^{(i+1)}$.
\\
 The use of Random Projection ANN, or ELMs, has a vast literature and includes theoretical foundings and applications to various fields; we refer to the review papers \cite{ding2015extreme,huang2015trends,wang2022review} for interested readers.\B

\textbf{Remark:}
The suggested ELM training procedure, based on a choice made via random selection for part of the unknown parameters and some linear problems to be solved for the determination of the remaining unknown parameters, has been very vastly studied in the literature.
When ELM was defined by Huang and co-authors, it gained huge popularity and lead to very interesting developments \cite{huang2006extreme}.
Nevertheless, the debate on the original denomination or on the first researcher that has introduced such a procedure populates both scientific papers and internet blogs.
In the authors' opinion, the contribution due to Huang is out of doubt, though it is far from the scope of the authors to crown one scholar or another.
We refer to two websites\footnote{ https://elmorigin.wixsite.com/originofelm ; https://www.extreme-learning-machines.org/} for those who are interested in such debate.

\textbf{Remark:} {
In the ELM framework, the multilayer construction (i.e., deep ANNs) becomes a linear combination of trial functions, which are, in turn, compositions of activation functions, themselves generated by the random choices of the parameters.
The deep architecture is usually adopted when considering high-dimensional problems.
However, in our experience, for the low-dimensional approximation problem considered, the structure of the deep network is not necessary to meet the required properties, but is only useful to change the shape of the functions.
}

\subsection{Resulting interpolation problem}
\label{subsec:Numerical problem}

Based on the ELM training procedure adopted in the present work, we can define a function $\tilde{u}$ to summarize the action of the entire network (called the \emph{network function}).
The network function relates the network input and the external weights with the network output, such that: $x^{(2)}=\tilde{u}(x^{(0)},\ul{w})$.
In our approximation problem, $x^{(0)}$ and $x^{(2)}$ are the independent and dependent variables of the target univariate scalar function sought to be interpolated.
The leading philosophy of the ELM framework is to learn only external weights  
since the following theorem holds see \cite{huang2006extreme}: 
\begin{thm}\label{Th3.1}
    Let $(x_j,y_j)\,,\ j=1,\dots, M $ be a set of data points such that $x_j\not=x_{j+1}$, and lets take the single-layer FFN with $N^{(1)}\le M$ neurons, then:
    \begin{equation}
        \tilde{u}(x_j; w_i)= \sum_{i=1}^{N^{(1)}} w_i  \psi_i^{(1)}(z_i^{(1)})
        \qquad
        \textnormal{with}
        \quad
        z_i^{(1)}=\kappa_i^{(1)}(x^{(0)},A_i^{(1)},\beta_i^{(1)})       
        \, ,
        \label{eq:ELM_net_function}
    \end{equation}
    Moreover, consider weights $A_i^{(1)}$ and bias $\beta^{(1)}$ randomly generated, independently from the data points $(x_j,y_j)$, according to any continuous probability distribution $\mathbb{P}$. Then, $\forall \varepsilon >0$, there exists a choice of $\ul{w}$ such that $\mathbb{P}(\| (\tilde{u}(x_j; w_i) - y_j)\| <\varepsilon)=1 $.
    Moreover, if $N^{(1)}=M$ then $w_i$ can be found such that  $\mathbb{P}(\| (\tilde{u}(x_j; w_i) - y_j) \| =0)=1 $.
\end{thm}


Under the single-hidden layer architecture and ELM training framework adopted in the present work, the underlying approximation space is formed by the activation functions of the hidden layer $\psi_i^{(1)}$ (with $i=1\dots N^{(1)}$).
A common nomenclature to refer to the functions $\psi_i^{(1)}$ is \emph{trial} functions.
The trial functions are fixed because the internal parameters of the activation functions are given after the random projection.
Following what presented in \cite{calabro2021extreme}, in our numerical experiments we randomly select internal weights and biases, within a certain fixed range, so that the the centers of the activation functions are located in random points.

In order to select the external weights $\ul{w}$, we introduce a loss function to be optimized.
Since we focus on the problem of function interpolation, we call $f$ a generic target function, consider $M$ training points $(x_j,y_j)$ such that $f(x_j):=y_j\,,\ j=1,\dots , M $, with $M\le N^{(1)}$ and we look for the network function $\tilde{u}$ \eqref{eq:ELM_net_function} as a linear combination of trial functions and such that: 
\begin{equation}\label{eq:loss}
    \tilde{u} (x_j,\ul{w})=y_j \ , \quad \forall j=1,\dots,M .
\end{equation}

The presence of free parameters acting only linearly clearly distinguishes the ELM case from other training approaches, in which the determination of the degrees of freedom is a nonlinear problem.
Accordingly, in networks trained via ELM, the exact (optimal) choice of weights $\ul{w}$ can be determined in a ``one--shot" process.
Notice that such a training can be also considered as an useful procedure for the initialization of the iterative optimization procedure of standard ANNs.

In fact, our target problem \eqref{eq:loss} can be rewritten in matrix form as follows.
\begin{equation}\label{eq:exactness}
  \mathbb{S} \cdot \ul{w} =\ul{y}\ ,
  \qquad
  \textnormal{with }
  \mathbb{S} = \begin{bmatrix} 
    \psi_1^{(1)}(x_1) & \dots  & \psi_{N^{(1)}}^{(1)}(x_1)\\
    \vdots & \ddots & \vdots\\
    \psi_1^{(1)}(x_M) & \dots  & \psi_{N^{(1)}}^{(1)}(x_M) 
    \end{bmatrix}
    , 
\end{equation}
where $\mathbb{S} \in \mathbb{R}^{M \times N^{(1)}}$ can be interpreted as a collocation matrix, since it consists of different functions (represented by each column of the matrix) evaluated in a set of points (each one represented by a row of the matrix).

Since Theorem \ref{Th3.1} guarantees (at least with probability 1) the existence of an interpolation network, the optimal solution of problem \eqref{eq:loss} is computed by the resolution of \eqref{eq:exactness}.
Then, we aim to study the behavior of the constructed interpolant in terms of function approximation as it is well known that interpolation on nodes does not imply (in general) good behavior between nodes see, e.g., \cite{trefethen2019approximation}.

Theorem \ref{Th3.1} states that square problems are full rank with probability 1.
Ill-conditioning can occur because of two reasons: the location of interpolating nodes or the presence of similar trial functions see, e.g., \cite{corless2020runge,trefethen2019approximation}.
In these cases, the rank can be (at least to machine precision) lowered.
A convenient way to overcome both these issues is to adopt an overparametrized system, i.e. the number of network parameters $N^{(1)}$ is higher than the number of data points $M$.  We point out that the requested interpolation conditions are  matched also in these cases. Then one of the solutions has to be selected because the linear system with $M<N^{(1)}$ is \emph{underdetermined}. We chose to solve it by the least square procedure that gives the minimum-norm solution, as is done by backslash command in MATLAB. This leads to an increased computational cost with respect to the square case, related to the construction  of the normal equations. When overparametrizing, we take $M=N^{(1)}/2$, following \cite{calabro2021extreme}. 
Nevertheless, in the next section where Runge's example is considered, we consider also the square case (i.e. $M=N^{(1)}$).

\section{Numerical results}
\label{sec:NumericalResults}
In this section, we study the problem of approximating a regular function through a single-layer ANN{ trained via ELM;  in particular, we study the obtained convergence properties increasing the number of training points $M$ (i.e. the interpolation nodes) and correspondingly  increasing the number of neurons $N^{(\mathcal{L})}$ (i.e. the number of degrees of freedom). In particular, we consider the square case, where $M=N^{(\mathcal{L})}$ and the overparametrized case, where $M=N^{(\mathcal{L})}/2$.  }
Our focus is on the behavior of the interpolation independently of the location of nodes, and we present tests done on three different sets of interpolation nodes: equispaced including boundary points, Chebychev points, and randomly generated points.
We compare such convergence with what is obtained with the use of degree-increasing polynomial interpolation on the same nodes.
The polynomial interpolation on Chebychev nodes is chosen as reference target behavior, since it is, on one side, (quasi-) optimal, as widely discussed in \cite{trefethen2008gauss,trefethen2019approximation} and, on the other side, efficiently implemented in the CHEBFUN suite \cite{driscoll2014chebfun}.
The results reported regarding polynomial interpolation (Poly) are obtained by the CHEBFUN routine $chebfun.interp1(x,y)$, while the training of  ANN is done via ELM as discussed in the previous section.

When adopting ANNs, we construct the interpolation network $\tilde{u}$ of \eqref{eq:loss} with various activation functions: summarizing what was presented in the previous section, the explored options are the following:
\begin{itemize}
    \item ANN-based interpolation with  Logistic Sigmoidal activation function (LS);
    \item ANN-based interpolation with Gaussian Radial Basis activation function (GRB);
    \item ANN-based interpolation with SoftPlus activation function (SP).
\end{itemize}

In all cases, we are interested in the behavior outside the interpolation nodes, where the exactness conditions are guaranteed by Theorem \ref{Th3.1}. Then, we introduce the interpolation error:
\begin{equation*}
    Err=\| \tilde{u}-f  \|_2
\end{equation*}
where the interpolating $\tilde{u}$ and the interpolated -target- $f$ are functions evaluated on 4000 equispaced points and the 2-norm ($\|\dots\|_2$) is evaluated via a trapezoidal rule.

We use different target functions with different regularities so that results are reported in proper scales according to the expected convergence: when available, we report also the reference convergence ratio. 

In Section \ref{subsec:Runge} we discuss results obtained in the case of Runge's example. For this function, we present various convergence plots. 
First of all, we consider the case where the ANN architecture is squared: the number of neurons in the hidden layer $N^{(\mathcal{L})}$ is the same as the number of interpolation nodes $M$. Then, we consider the case where the ANN architecture is overparametrized: in this case, we choose $M=N^{(\mathcal{L})}/2$.  In the square case, we notice a good but not optimal convergence, while in the latter we observe a convergence that exactly meets the reference one, i.e. the degree-increasing polynomial interpolation on Chebychev nodes. We also present convergence results regarding the approximation of the derivative: we train the ANN on the nodes and evaluate the exact derivative of the trained ANN, then we estimate the error with respect to the exact derivative of the interpolated function. Also in this case we obtain very good results. 

In Section \ref{subsec:OtherFunctions} we discuss results obtained in other test cases, focusing on overparametrized ANNs. We range our tests following what is presented in \cite{trefethen2019approximation}: both analytic functions and functions with lower regularities with known behavior of the error  are considered. In all cases, good results are obtained, especially due to the fact that the convergence is observed independently of the location of the interpolating nodes.
 This result demonstrates that the interpolation is stable with respect to the choice of nodes. In the current study, we have focused our attention on exact interpolation, thus not considering noisy functions. These cases could be better handled by adding regularization terms; however, it is important to note that the inclusion of regularization terms does not allow for exact interpolation. \B

\subsection{Runge's function approximation}
\label{subsec:Runge}
As a first test, we consider a 
well-known benchmark function based on the so-called Runge's example. The target problem is to approximate the following bell-shaped function:
\begin{equation}
    f_R(x)=\frac{1}{1+25x^2}, \qquad x \in [-1,1].
    \label{eq:runge}
\end{equation}
It is well-known that for this function polynomial interpolation on equispaced points diverges while adopting an interpolation based  on Chebychev nodes results in a geometrical convergence  with respect to the number of interpolating points $M$, that is, the convergence rate is $O(C^{-M})$ for some constant $C > 1$ \cite{platte2011impossibility,10.1093/imanum/dry024,boyd2009exponentially}; in particular, calling $p^{(C)}_M$ the interpolating polynomial on $M$ Chebychev nodes, then: \begin{equation}\label{eq:errChebRunge}
    \| p^{(C)}_M-f_R  \|=o\left( \left( \dfrac{1+\sqrt{26}}{5} \right)^{-M}\right)\ .
\end{equation}


\begin{figure}[ht]%
\centering
\includegraphics[width=0.335\textwidth  ,keepaspectratio]{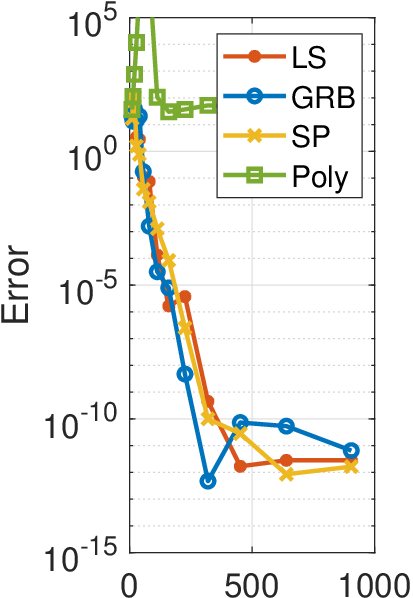}
\includegraphics[width=0.32\textwidth,keepaspectratio]{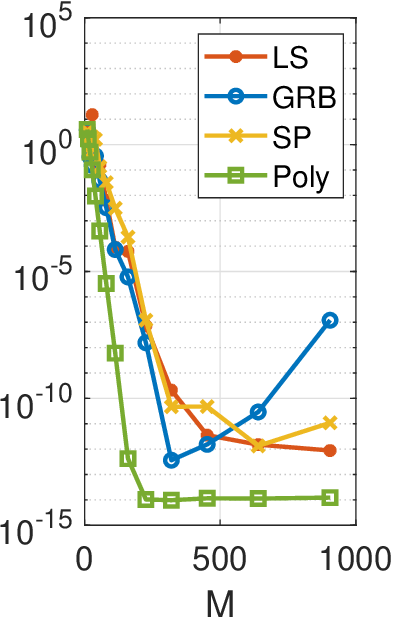}
\includegraphics[width=0.32\textwidth,keepaspectratio]{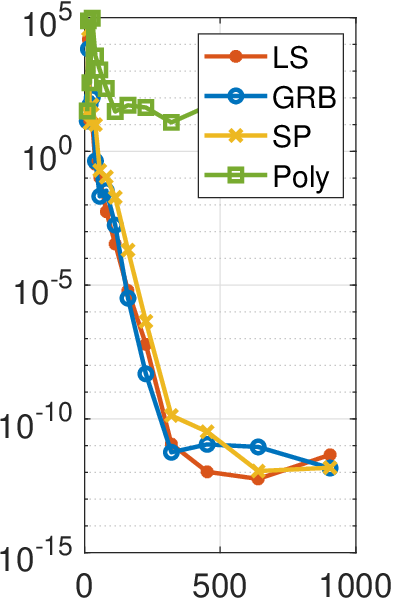}
\caption{Runge's example: computed error, square case $M=N^{(\mathcal{L})}$,  where $M$ is the number of nodes and $N^{(\mathcal{L})}$ is the number of neurons in the hidden layer. The left panel is the case of equispaced nodes, the center is the case of Chebychev nodes, and the right is the case of randomly generated nodes. Convergence is tested at different choices of the activation function and increasing $M$ \B and compared with one of the polynomial interpolation. In this case, the linear problem \eqref{eq:exactness} is squared. \label{fig1}}
\end{figure}
\begin{table}[] 
\begin{tabular}{lllll}
M      & LS                & RB                & SP                & Poly              \\
\_\_\_ & \_\_\_\_\_\_\_\_\_\_ & \_\_\_\_\_\_\_\_\_\_ & \_\_\_\_\_\_\_\_\_\_ & \_\_\_\_\_\_\_\_\_\_ \\
10     & 5.1700e+01           & 2.0200e+01           & 1.6600e+01           & 3.6700e+01           \\
20     & 7.0200e+01           & 3.8500e+01           & 1.1300e+02           & 7.5400e+02           \\
40     & 2.7900e+00           & 6.2100e-01           & 7.7700e+01           & 8.8200e+05           \\
80     & 7.0300e-03           & 9.1900e-03           & 1.3200e-02           & 9.2000e+06           \\
160    & 3.3800e-05           & 9.2600e-06           & 4.4300e-05           & 7.0000e+01           \\
320    & 6.2200e-10           & 1.5500e-13           & 3.5500e-11           & 3.8800e+01          
\end{tabular}
\caption{  Computed errors of the interpolation of Runge's function increasing the number of nodes $M$ - in the first column -  both for different choices of the activation function and polynomial interpolation. Notation follows what presented at the beginning of Section 3. The table reports the errors computed for equispaced nodes in the square case $M=N^{(\mathcal{L})}$, where $N^{(\mathcal{L})}$ is the number of neurons in the hidden layer. Such results are those depicted in the first panel of Figure \ref{fig1}. }
\label{Tab:1}
\end{table}

\begin{table}[] 
\begin{tabular}{lllll}
M      & LS                & RB                & SP                & Poly              \\
\_\_\_ & \_\_\_\_\_\_\_\_\_\_ & \_\_\_\_\_\_\_\_\_\_ & \_\_\_\_\_\_\_\_\_\_ & \_\_\_\_\_\_\_\_\_\_ \\
10     & 3.6300e+00           & 3.8000e+00           & 3.4000e+00           & 3.9200e+00           \\
20     & 3.5200e-01           & 3.3000e-01           & 4.4800e-01           & 5.0600e-01           \\
40     & 8.8000e+00           & 9.7900e-02           & 4.0400e-01           & 9.5000e-03           \\
80     & 1.0400e-02           & 2.9000e-02           & 3.8100e-02           & 3.3600e-06           \\
160    & 3.1300e-06           & 8.9400e-07           & 2.1100e-04           & 4.1900e-13           \\
320    & 5.1000e-11           & 1.2400e-12           & 4.5500e-10           & 1.6500e-14          
\end{tabular}
\caption{ Computed errors for the Runge's function. For the notation, see Table \ref{Tab:1}. In this case, we use Chebychev nodes, as reported in the second panel of Figure \ref{fig1}.}
\label{Tab:2}
\end{table}


\begin{table}[] 
\begin{tabular}{lllll}
M      & LS                & RB                & SP                & Poly              \\
\_\_\_ & \_\_\_\_\_\_\_\_\_\_ & \_\_\_\_\_\_\_\_\_\_ & \_\_\_\_\_\_\_\_\_\_ & \_\_\_\_\_\_\_\_\_\_ \\
10     & 1.5500e+01           & 1.3500e+01           & 2.9200e+01           & 3.2100e+01           \\
20     & 2.6900e+01           & 5.5900e+01           & 3.0200e+01           & 3.6600e+02           \\
40     & 1.4200e+00           & 1.5500e+00           & 4.2800e-01           & 3.6500e+03           \\
80     & 9.8600e-03           & 3.9100e-02           & 2.3600e-02           & 3.0200e+02           \\
160    & 7.1700e-05           & 4.4300e-06           & 1.9600e-04           & 4.0400e+01           \\
320    & 5.8900e-12           & 2.3500e-12           & 8.5400e-11           & 1.1800e+01          
\end{tabular}
\caption{ Computed errors for the Runge's function. For the notation, see Table \ref{Tab:1}. In this case, we use random nodes, as reported in the third panel of Figure \ref{fig1}. }
\label{Tab:3}
\end{table}

\begin{figure}[ht]%
\centering
\includegraphics[width=0.335\textwidth  ,keepaspectratio]{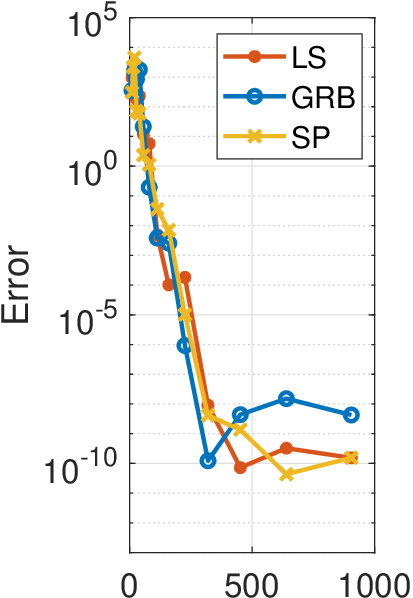}
\includegraphics[width=0.32\textwidth  ,keepaspectratio]{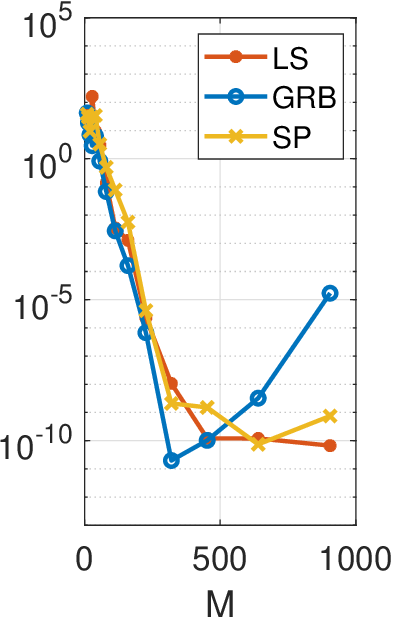}
\includegraphics[width=0.32\textwidth  ,keepaspectratio]{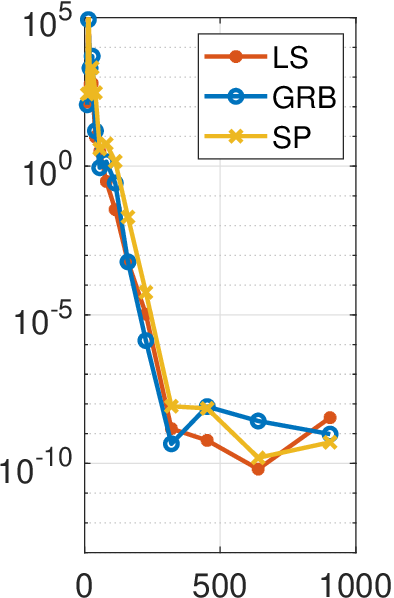}
\caption{Runge's example: computed error for the derivative, square case $M=N_{\mathcal{L}}$  where $M$ is the number of nodes and $N^{(\mathcal{L})}$ is the number of neurons in the hidden layer\B. The left panel is the case of equispaced nodes, the central panel is the case of Chebychev nodes, and the right panel is the case of randomly generated nodes at different choices of the activation function. The difference is between the exact derivative of a function $f_R$ in \eqref{eq:runge} and the computed derivative of the function obtained by interpolation, the one used in Figure \ref{fig1}. \label{fig2}}
\end{figure}

\begin{figure}[ht]%
\centering
\includegraphics[width=0.335\textwidth  ,keepaspectratio]{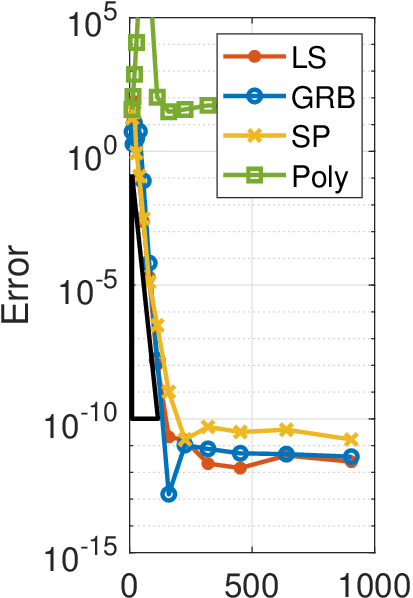}
\includegraphics[width=0.32\textwidth  ,keepaspectratio]{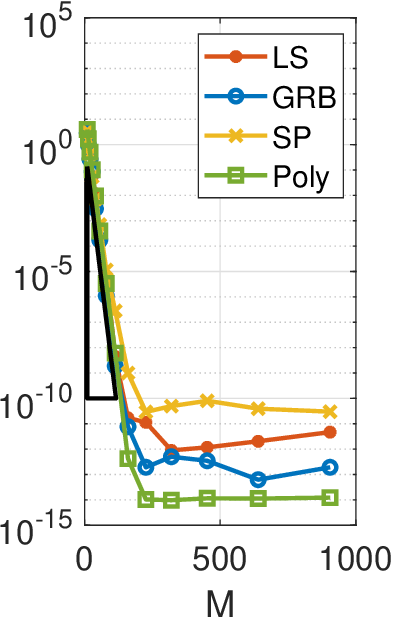}
\includegraphics[width=0.32\textwidth  ,keepaspectratio]{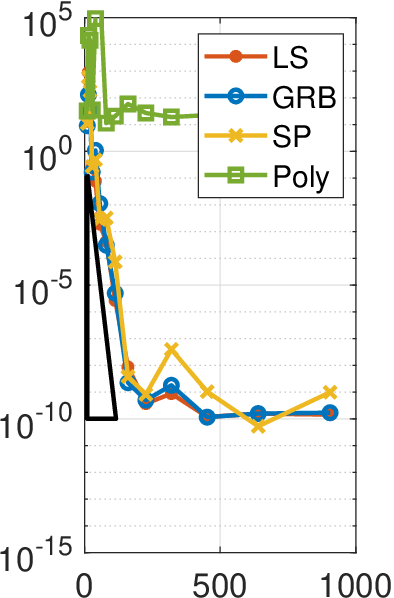}
\caption{Runge's example: computed error with respect to the number of interpolating nodes $M$, overparametrized case with $M=N^{(\mathcal{L})}/2$  where $N^{(\mathcal{L})}$ is the number of neurons in the hidden layer\B. The left panel is the case of equispaced nodes, the central panel is the case of Chebychev nodes, the right panel is the case of randomly generated nodes. The linear problem \eqref{eq:exactness} is solved by least squares. The black triangle is the reference convergence  (interpolating polynomial on Chebychev nodes), as reported in \eqref{eq:errChebRunge}. \label{fig3}}
\end{figure}

\begin{figure}[ht]%
\centering
\includegraphics[width=0.335\textwidth  ,keepaspectratio]{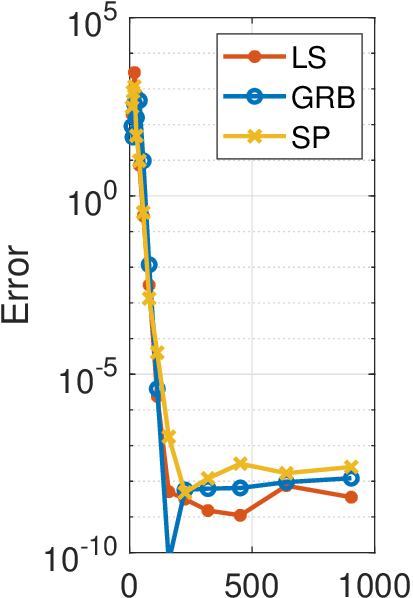}
\includegraphics[width=0.32\textwidth  ,keepaspectratio]{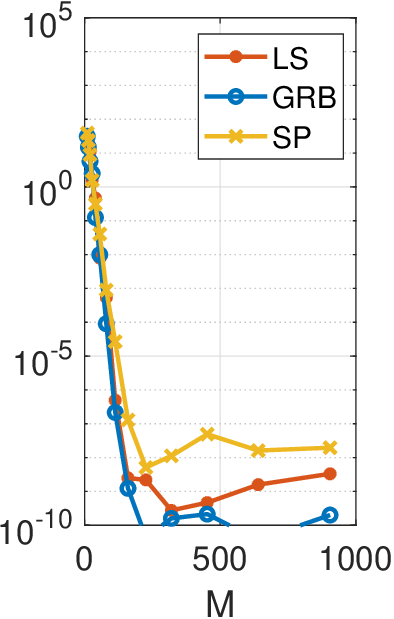}
\includegraphics[width=0.32\textwidth  ,keepaspectratio]{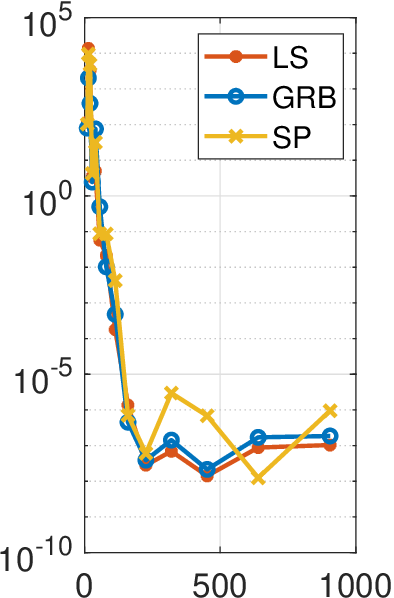}
\caption{Runge's example: computed error for the derivative with respect to the number of interpolating nodes $M$, overparametrized case with $M=N^{(\mathcal{L})}/2$  where  $N^{(\mathcal{L})}$ is the number of neurons in the hidden layer\B. The left panel is the case of equispaced nodes, the central panel is the case of Chebychev nodes, and the right panel is the case of randomly generated nodes. The difference is between the exact derivative of function $f_R$ in \eqref{eq:runge} and the computed derivative of the function obtained by interpolation, the one used in Figure \ref{fig3}. \label{fig4}}
\end{figure}

In figures \ref{fig1}-\ref{fig4} we report the errors obtained when testing our procedure with the target function being Runge's one and its derivative. In the first two, figures \ref{fig1}-\ref{fig2} the  convergence is done in the squared case $M=N^{(\mathcal{L})}$, respectively interpolating the function \eqref{eq:runge} and its derivative.
In the first case, we can notice that the convergence is exponential although not optimal, i.e. geometrical but with a different constant.
Surprisingly, exponential convergence is independent of the choice of the interpolation nodes.
Convergence in all cases is obtained up to a precision that is related to the condition of the linear problem \eqref{eq:exactness} that is solved.
We point out that no particular care has been taken for the resolution of such a problem: all tests are done via the backslash command in MATLAB: recall that the determination of the parameters of the ANN does not involve optimization: the internal parameters are fixed via ELM and the external weights are solved through a square linear problem.  

In the second case, the one reported in Figure \ref{fig2}, we plot the errors obtained for the approximation of the derivative of Runge's function in the squared case. We compute the interpolant as in the first case and then evaluate the error that is done with respect to the derivative of the original functions. This can be done in our case because both the exact derivative of the function $f_R$ and the exact derivative of the interpolating ANN can be done in an exact way. From the reported test we can conclude that the convergence of the network function is exponential also in this case, which reveals a surprising ability of such interpolating ANN to approximate very well also the derivative of the function.  Polynomial interpolation is not well suited for the approximation of the derivative of an interpolated function so the results of this test are not reported in the figures.

In order to recover the exact optimal convergence for Runge's example, we consider the overparametrized ANN and fix $M=N^{(\mathcal{L})}/2$.
Recall that the interpolation problem \eqref{eq:exactness} for the single-layer ANN \eqref{eq:ELM_net_function} remains linear in the unknowns - free parameters - $w$ so that the training is reduced to the resolution of an underdeterminated system, ad ex. via least squares. 
In Figures \ref{fig3} and \ref{fig4} we report the results obtained in this overparametrized case for the function itself and for its derivative, respectively.
Error convergence is optimal in this case, see equation \eqref{eq:errChebRunge}, and we notice again that is independent of the choice of the interpolation nodes. To emphasize this, we report in a black triangle the reference
convergence, interpolating polynomials on Chebychev nodes. The least-square procedure, then, adapts the solution in this overparametrized setting, so that the enrichment of the space is effectively used even if the number of interpolating points does not match the number of independent functions used for approximation.
For these reasons, we propose this overparametrization as the procedure for all the next numerical tests.

\begin{table}[] 
\begin{tabular}{lllll}
M      & LS                & RB                & SP                & Poly              \\
\_\_\_ & \_\_\_\_\_\_\_\_\_\_ & \_\_\_\_\_\_\_\_\_\_ & \_\_\_\_\_\_\_\_\_\_ & \_\_\_\_\_\_\_\_\_\_ \\
10     & 1.0500e+01           & 5.3900e+00           & 1.6900e+01           & 3.6700e+01           \\
20     & 7.2200e+01           & 1.0300e+01           & 2.9000e+01           & 7.5400e+02           \\
40     & 8.7200e-02           & 5.6200e+00           & 1.1900e-01           & 8.8200e+05           \\
80     & 2.0300e-05           & 6.6600e-05           & 1.2800e-05           & 9.2000e+06           \\
160    & 2.1000e-11           & 3.1300e-13           & 1.0200e-09           & 7.0000e+01           \\
320    & 2.1300e-12           & 7.5300e-12           & 4.9900e-11           & 3.8800e+01          
\end{tabular}
\caption{ Computed errors of the interpolation of Runge's function increasing the number of nodes $M$ - in the first column -  both for different choices of the activation function and polynomial interpolation. Notation follows what presented at the beginning of Section 3. The table reports the errors computed for equispaced nodes in the overparametrized case $M=N^{(\mathcal{L})}/2$, where $N^{(\mathcal{L})}$ is the number of neurons in the hidden layer. Such results are those depicted in the first panel of Figure \ref{fig3}.}
\label{Tab:4}
\end{table}

\begin{table}[]
\begin{tabular}{lllll}
M      & LS                & RB                & SP                & Poly              \\
\_\_\_ & \_\_\_\_\_\_\_\_\_\_ & \_\_\_\_\_\_\_\_\_\_ & \_\_\_\_\_\_\_\_\_\_ & \_\_\_\_\_\_\_\_\_\_ \\
10     & 3.0300e+00           & 2.8300e+00           & 3.4900e+00           & 3.9200e+00           \\
20     & 5.1800e-01           & 2.6400e-01           & 4.2000e-01           & 5.0600e-01           \\
40     & 1.1400e-02           & 3.0000e-03           & 7.0400e-03           & 9.5000e-03           \\
80     & 6.6300e-06           & 1.1100e-06           & 1.1500e-05           & 3.3600e-06           \\
160    & 1.7000e-11           & 7.6500e-12           & 9.9500e-10           & 4.1900e-13           \\
320    & 8.9100e-13           & 5.0000e-13           & 4.9500e-11           & 1.6500e-14          
\end{tabular}
\caption{ Computed errors for the Runge's function. For the notation, see Table \ref{Tab:4}. In this case, we use Chebychev nodes, as reported in the second panel of Figure \ref{fig3}.}
\label{Tab:5}
\end{table}

\begin{table}[]
\begin{tabular}{lllll}
M      & LS                & RB                & SP                & Poly              \\
\_\_\_ & \_\_\_\_\_\_\_\_\_\_ & \_\_\_\_\_\_\_\_\_\_ & \_\_\_\_\_\_\_\_\_\_ & \_\_\_\_\_\_\_\_\_\_ \\
10     & 3.0300e+00           & 2.8300e+00           & 3.4900e+00           & 3.9200e+00           \\
20     & 5.1800e-01           & 2.6400e-01           & 4.2000e-01           & 5.0600e-01           \\
40     & 1.1400e-02           & 3.0000e-03           & 7.0400e-03           & 9.5000e-03           \\
80     & 6.6300e-06           & 1.1100e-06           & 1.1500e-05           & 3.3600e-06           \\
160    & 1.7000e-11           & 7.6500e-12           & 9.9500e-10           & 4.1900e-13           \\
320    & 8.9100e-13           & 5.0000e-13           & 4.9500e-11           & 1.6500e-14          
\end{tabular}

\caption{ Computed errors for the Runge's function. For the notation, see Table \ref{Tab:4}. In this case, we use random nodes, as reported in the third panel of Figure \ref{fig3}. }
\label{Tab:6}
\end{table}

\B
\clearpage

\subsection{Other function approximation tests}
\label{subsec:OtherFunctions}

In this section, we present some examples to test our procedure with the target function being various regular functions.

Notation follows what was presented in the previous section.
The reference interval is always $[-1,1]$.
The considered test covers all possible behaviors of regular functions and we consider: \begin{itemize} \item a polynomial function (Section \ref{Sub1}); \item analytic functions of the three possible types: analytic in the whole complex plane, analytic in the reference domain with a pole far away from the domain, analytic with poles in the complex plane near the considered domain (Section \ref{Sub2}); \item functions with lower regularities (Section \ref{Sub3}). \end{itemize} 
Moreover, while selecting such functions, we include also highly oscillating and steep gradient behaviors: all examples in sections \ref{Sub2}-\ref{Sub3} follow the work reported in \cite{trefethen2019approximation}.

\subsubsection{Polynomial function}
\label{Sub1}

In this section, we present a polynomial  example to test the accuracy and stability of our procedure. We generate, with a random procedure based on the $randi$ function of Matlab, the following:\begin{multline}
 f(x)= 4x^{16}-4x^{15}+ 8x^{13} +8x^{12} -8x^{11} -6x^{10} -9x^{9}-x^{8}-10x^{7}+ \\ -x^{6}+ 3x^{5}-5x^{4} +4x^{3}+ 2x^{2} -10x+ 1 
\end{multline}

\begin{figure}[ht]%
\centering
\includegraphics[width=0.335\textwidth  ,keepaspectratio]{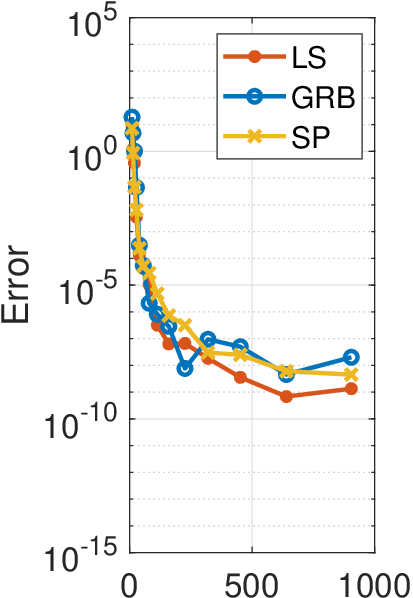}
\includegraphics[width=0.32\textwidth  ,keepaspectratio]{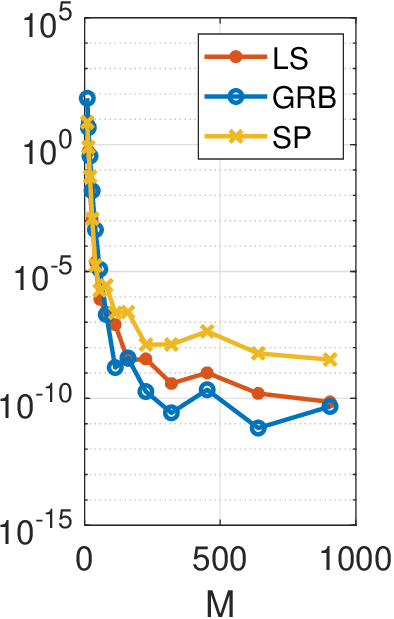}
\includegraphics[width=0.32\textwidth  ,keepaspectratio]{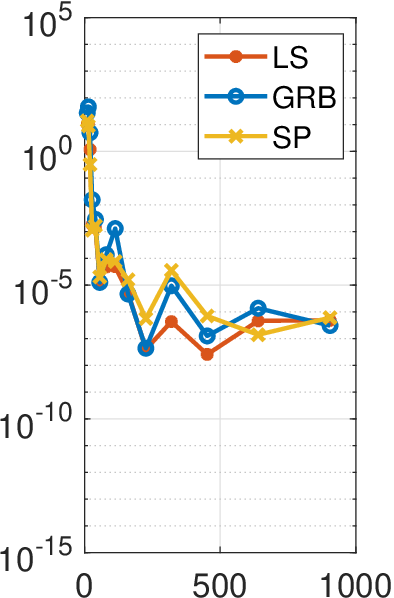}
\caption{Interpolation of a polynomial function: computed error, overparametrized case with $M=N^{(\mathcal{L})}/2$  where $M$ is the number of nodes and $N^{(\mathcal{L})}$ is the number of neurons in the hidden layer\B. Convergence is tested at different choices of the activation function and increasing 
$M$. The left panel is the case of equispaced nodes, the central panel is the case of Chebychev nodes, and the right panel is the case of randomly generated nodes. The linear problem \eqref{eq:exactness} is solved by least squares. Convergence is reported in a semi-log scale.  \label{fig_polrandom}}
\end{figure}

%
%
%
%
%
%
%
%
%
This test is of particular interest because our choice of network does not reproduce polynomials, thus, as expected, also in this case we can evaluate convergence but not exactness. Then, the aim is to test the accuracy of our procedure and, by considering increasing numbers of degrees of freedom, the stability. In Figure \ref{fig_polrandom} we notice that the convergence is exponential, and stops at different accuracy levels, which is expected due to the difference in condition number when the location of nodes is changed. Also, stability is confirmed, being in all tests the higher accuracy is maintained when the number of nodes increases. As already noticed in Runge's example, we obtain good accuracy in all tests, and independently of the location of nodes as we present tests done on three different sets of
interpolation nodes and the convergence is maintained.

\clearpage

\subsubsection{Other analytic functions}
\label{Sub2} 
In this section, we present some examples to test the accuracy and stability of our procedure when applied to target functions that are analytic functions in the reference interval.
 In particular, the functions selected are those investigated in \cite{trefethen2019approximation}, in order to have benchmark results on the convergence rate and we compare the convergence of our network function with the geometric convergence of the polynomial interpolation on Chebychev nodes.

\begin{figure}[ht]%
\centering
\includegraphics[width=0.335\textwidth  ,keepaspectratio]{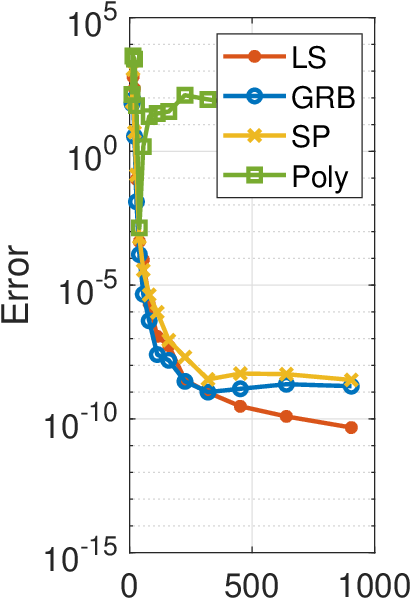}
\includegraphics[width=0.32\textwidth  ,keepaspectratio]{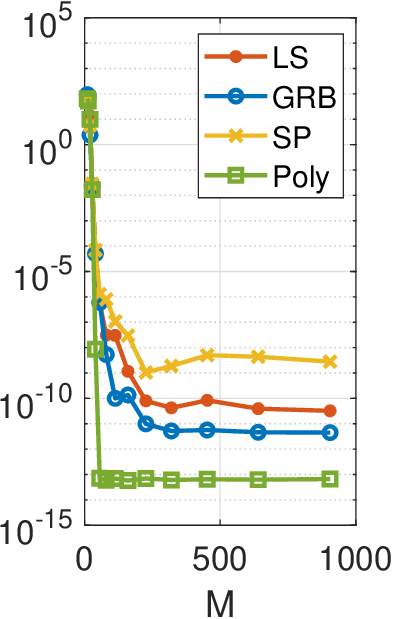}
\includegraphics[width=0.32\textwidth  ,keepaspectratio]{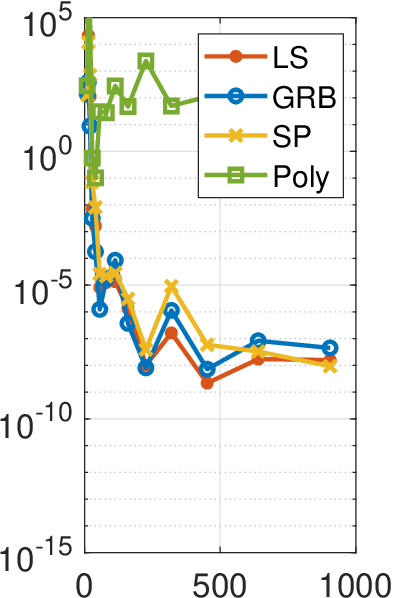}
\caption{Interpolation of the rapidly oscillating function $cos(20 x)$: computed approximation error, overparametrized case with $M=N^{(\mathcal{L})}/2$  where $M$ is the number of nodes and $N^{(\mathcal{L})}$ is the number of neurons in the hidden layer\B. Convergence is tested at different choices of the activation function and increasing 
$M$ and compared with one of the polynomial interpolation. The expected convergence in the reference case is more than exponential. The left panel is the case of equispaced nodes, the central panel is the case of Chebychev nodes, and the right panel is the case of randomly generated nodes. The linear problem \eqref{eq:exactness} is solved by least squares. Convergence is reported in a semi-log scale.   \label{fig_coseno}}
\end{figure}

The first example that we present is $f(x) = \cos(20 x)$ that is analytic not just on $[-1, 1]$, but in the whole complex plane, and the convergence expected for the polynomial interpolation on the Chebychev nodes is even faster than geometric. Moreover, this function is rapidly oscillating in $[-1, 1]$. In Figure \ref{fig_coseno} we report the computed error.  We observe again that by applying our procedure we  obtain convergence independently on the choice of the interpolation nodes. Also, as before, we notice that the convergence stops at different accuracy levels, due to the condition number of location of the interpolation nodes.

\begin{figure}[ht]%
\centering
\includegraphics[width=0.335\textwidth  ,keepaspectratio]{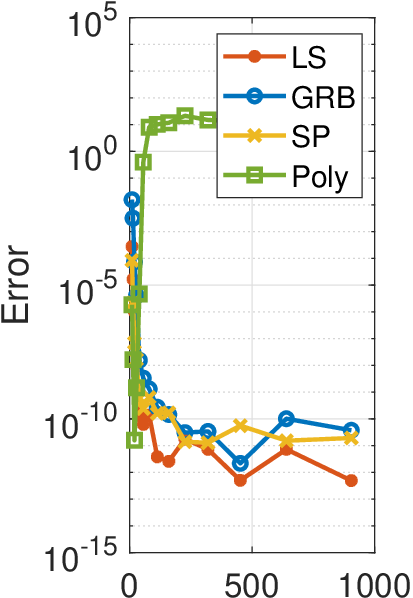}
\includegraphics[width=0.32\textwidth  ,keepaspectratio]{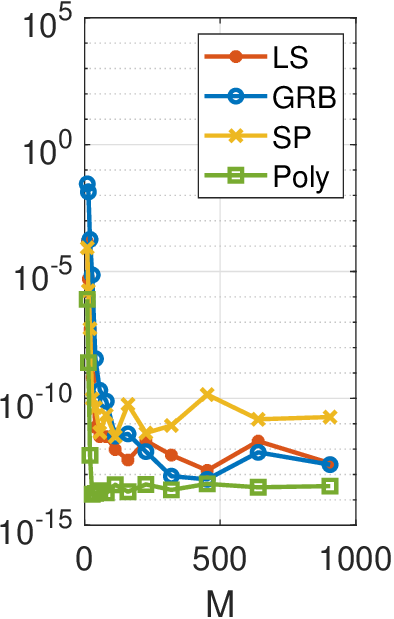}
\includegraphics[width=0.32\textwidth  ,keepaspectratio]{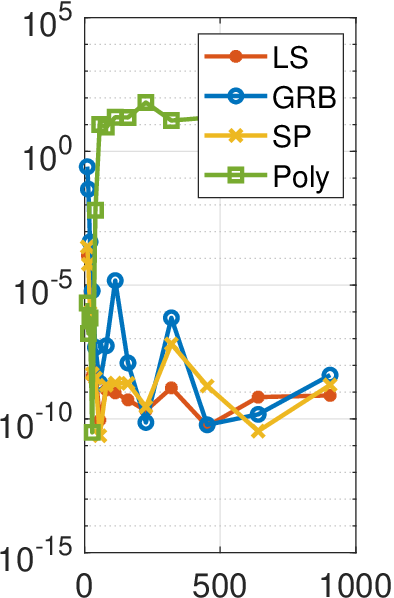}
\caption{Interpolation of the analytic function with a real singularity $\sqrt{(2-x)}$: computed error with respect to the number of interpolating nodes $M$, overparametrized case with $M=N^{(\mathcal{L})}/2$  where  $N^{(\mathcal{L})}$ is the number of neurons in the hidden layer\B. The left panel is the case of equispaced nodes, the central panel is the case of Chebychev nodes, and the right panel is the case of randomly generated nodes. The linear problem \eqref{eq:exactness} is solved by least squares. Convergence is reported in a semi-log scale. \label{fig_radice}}
\end{figure}

\clearpage
The second example that we present is for the target function $f(x) = \sqrt{2-x}$, which is an example of a function with a real singularity, that has a branch outside of the considered interval, and relatively far away.
We obtain good accuracy in all tests reported in Figure \ref{fig_radice}.  As before, we can notice again that the procedure is stable and that the convergence is independent of the choice of the location of interpolation nodes and also that the convergence stops at different accuracy levels, due to the condition number of location of the interpolation nodes.

\begin{figure}[ht]%
\centering
\includegraphics[width=0.34\textwidth  ,keepaspectratio]{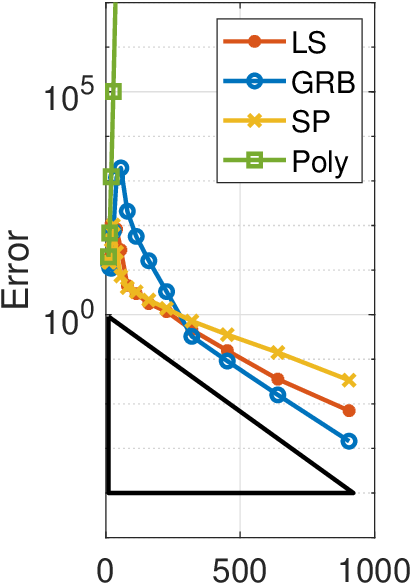}
\includegraphics[width=0.305\textwidth  ,keepaspectratio]{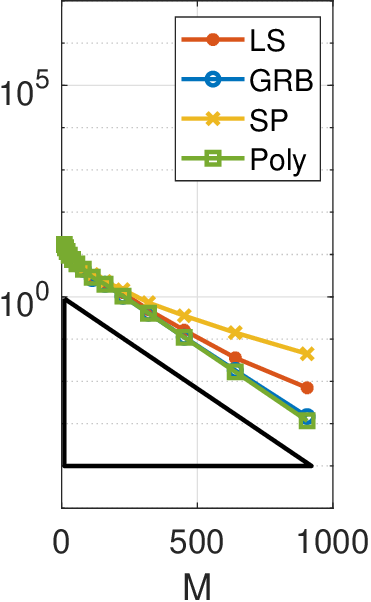}
\includegraphics[width=0.305\textwidth  ,keepaspectratio]{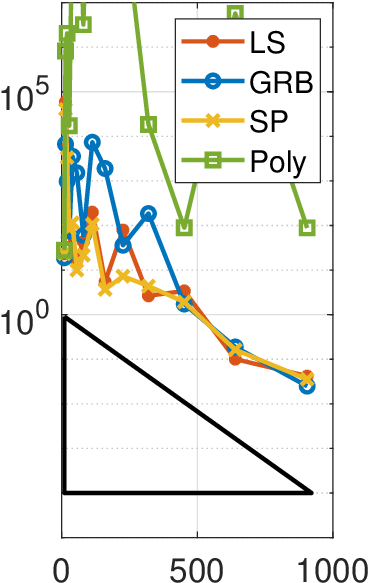}
\caption{Interpolation of the peaked function $\tanh(50\pi x)$: computed error with respect to the number of interpolating nodes $M$, overparametrized case with $M=N^{(\mathcal{L})}/2$  where $N^{(\mathcal{L})}$ is the number of neurons in the hidden layer\B. The left panel is the case of equispaced nodes, the central panel is the case of Chebychev nodes, and the right panel is the case of randomly generated nodes. The linear problem \eqref{eq:exactness} is solved by least squares. Convergence is reported in a semi-log scale. The black triangle in the reference convergence (interpolating polynomial on Chebychev nodes). \label{fig_tanh}}
\end{figure}

As a last example of an analytic function, we test the sharp gradient function $f(x) = \tanh(50\pi x)$.
This function is a more extreme but entirely analogous example of the Runge function.
Indeed it has two poles, but at $\pm 0.01i$, thus closer to the analyzed interval $[-1, 1]$. This makes the polynomial interpolation on Chebychev nodes, for the theory explained in \cite{trefethen2019approximation}, much slower, though still robust. Figure \ref{fig_tanh} shows also in this case that our procedure shows the reference convergence -reported in the black triangle- of the interpolating polynomial on Chebychev nodes, and this convergence rate is independent of the choice of the location of interpolation nodes.

\clearpage

\subsubsection{Differentiable functions}
\label{Sub3}
We now present some examples to test the accuracy and stability of our procedure under different regularity hypotheses of the target functions.
We use differentiable functions again selected among the functions examined in \cite{trefethen2019approximation} to study the smoothness-approximability link.
More precisely, the analyzed functions are, for an integer $\nu\geq 0$, $\nu-1$ times differentiable and with the $(\nu-1)$--th derivative absolutely continuous and $\nu$--th derivative of bounded variation in $[-1,1]$. From the results of Theorem 7.2 in \cite{trefethen2019approximation}, this implies the convergence at the algebraic rate $O(M^{-\nu})$ of the polynomial interpolation on Chebychev nodes. In these cases, we pass from an exponential to a polynomial convergence rate with respect to the number of nodes.\\ 

\begin{figure}[ht]%
\centering
\includegraphics[width=0.34\textwidth  ,keepaspectratio]{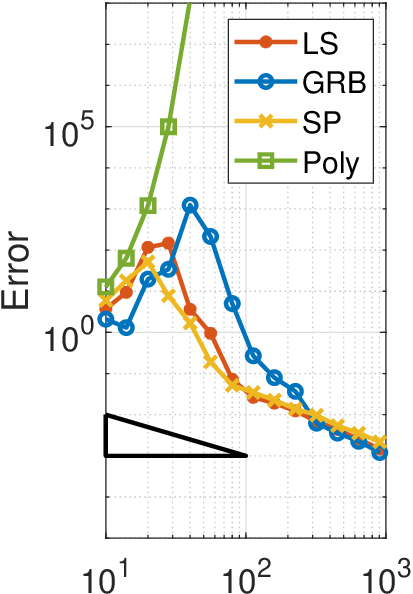}
\includegraphics[width=0.305\textwidth  ,keepaspectratio]{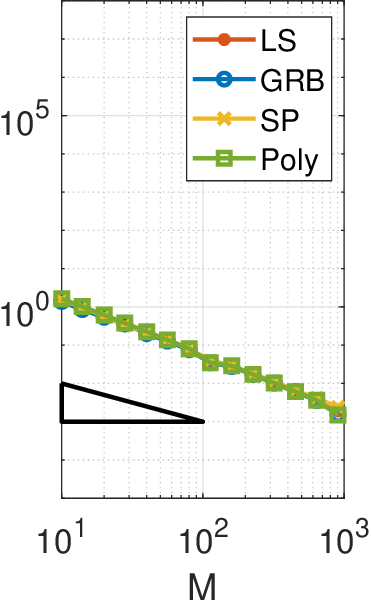}
\includegraphics[width=0.305\textwidth  ,keepaspectratio]{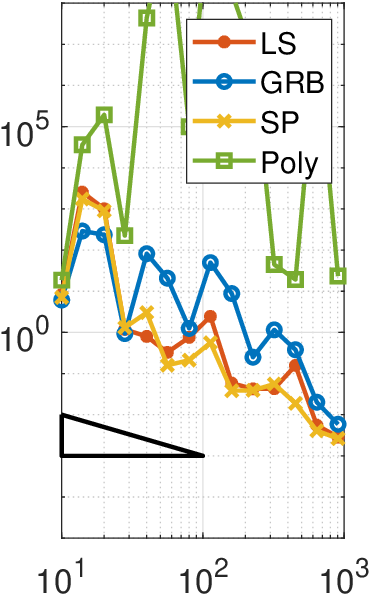}
\caption{Interpolation of the $C^0$ function $\left|{x}\right|$: computed error with respect to the number of interpolating nodes $M$, overparametrized case with $M=N^{(\mathcal{L})}/2$  where  $N^{(\mathcal{L})}$ is the number of neurons in the hidden layer\B. The left panel is the case of equispaced nodes, the central panel is the case of Chebychev nodes, and the right panel is the case of randomly generated nodes. The linear problem \eqref{eq:exactness} is solved by least squares. Convergence is reported in a log-log scale. The black triangle in the reference convergence (interpolating polynomial on Chebychev nodes). \label{fig_valass}}
\end{figure}

The first function that we use is $f(x) = \left|{x} \right|$, which is once differentiable with a jump in the first derivative at $x = 0$, thus it belongs to the space $C^0([-1,1])$. In Figure \ref{fig_valass} we plot the computed error in a log-log scale and the convergence curve of our procedure nicely matches $M^{-1}$ in all tests. Thus, not only it matches the theoretic convergence expected on Chebychev nodes, but also we find the same rate on equispaced and random nodes where the polynomial interpolation diverges. 

\begin{figure}[ht]%
\centering
\includegraphics[width=0.33\textwidth  ,keepaspectratio]{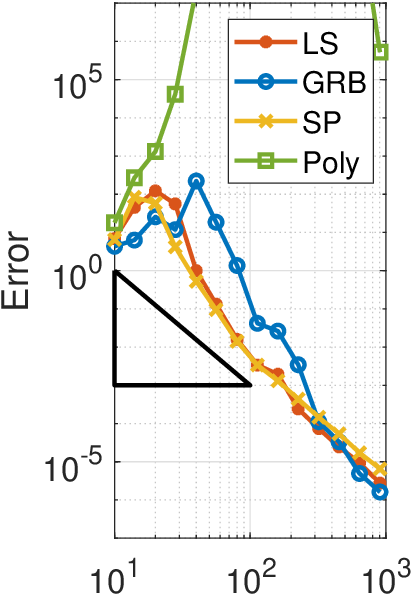}
\includegraphics[width=0.305\textwidth  ,keepaspectratio]{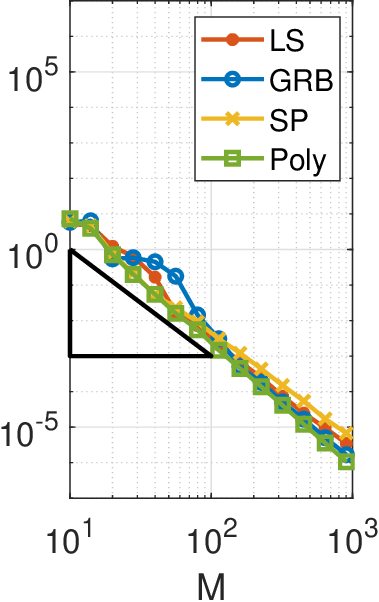}
\includegraphics[width=0.305\textwidth  ,keepaspectratio]{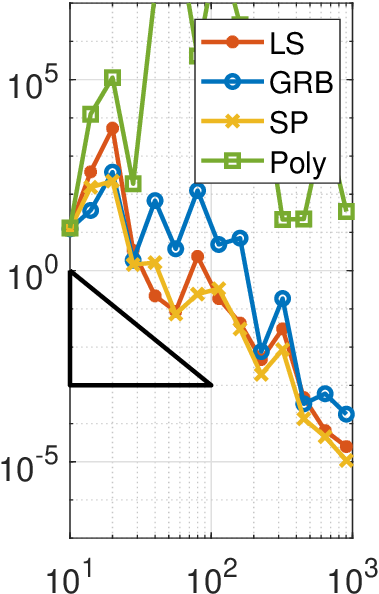}
\caption{Interpolation of the $C^2$ function $\left|\sin{5x}\right|^3$: computed error with respect to the number of interpolating nodes $M$, overparametrized case with $M=N^{(\mathcal{L})}/2$. The left panel is the case of equispaced nodes, the central panel is the case of Chebychev nodes, and the right panel is the case of randomly generated nodes. The linear problem \eqref{eq:exactness} is solved by least squares. Convergence is reported in a log-log scale. The black triangle is the reference convergence (interpolating polynomial on Chebychev nodes). \label{fig_c3}}
\end{figure}

The second function that we use is $f(x) = \left|\sin 5x \right|^3$ which is three times differentiable with jumps in the third derivative at $x = 0$ and $x= \pm \pi/5$, thus it belongs to the space $C^2([-1,1])$. In Figure \ref{fig_c3} we report the errors obtained in a log-log scale and we observe a similarity with what is predicted for the polynomial interpolation on Chebychev nodes, finding that the convergence curve of our procedure nicely matches $M^{-3}$ in all tests. Again, the same accuracy is observed for all the choices of the interpolation nodes.

\section{Conclusion}
\label{sec:Conclusion}

In the present paper, we have solved the problem of function approximation via interpolation within a class of single-layer artificial neural networks.
Such class of ANN is trained via a procedure, named ELM, that involves a random determination of internal parameters, given that the determination of the free parameters is a linear problem.
We conclude that the ANN interpolating function constructed as overparametrized ANN trained via ELM is a very good approximating function: the approximation error converges in most cases with the same rate with respect to the reference method when increasing the number of neurons, while keeping linear the approximation problem to be solved.
The first example that is taken into account in the present paper is Runge's function, where no divergence shows up also taking equispaced or randomly generating nodes.
In this and all other test cases, also rapidly oscillating and steep, the function is properly approximated and the approximation error decreases.
We suggest that if the target function is regular there is no need for deep structure in the ANN nor the training of internal parameters.

Our result, which shows that the Runge phenomenon is defected by ANN trained via ELM, opens up two important scenarios from our perspective. On the one hand, it highlights that the more general problem of overfitting in ANNs can be studied without the need to select training points. Instead, one can focus on the optimization procedure in overparametrized networks: the minimum norm solution without tuning internal parameters achieves the desired outcome. On the other hand, it allows for the more confident use of these functions in scientific machine learning contexts where regular solutions are expected. This provides a valid alternative to using Lagrange interpolation when global methods are desired. It allows for increasing the number of degrees of freedom without introducing spurious oscillations, while still achieving optimal convergence compared to global polynomial interpolation. 

Furthermore, as a future development, we believe that some of the results presented in this article can be revisited in light of the theoretical findings contained in \cite{neufeld2023universal} that we found as a recently published preprint while finalizing the first round of review of this work.

All included in the present paper adds informations on the use of ANN trained vial ELMs, that have been succesfully used in various applications, see \cite{ding2015extreme,huang2015trends,wang2022review}.
\B




\section*{Declarations}
\textbf{Compliance with Ethical Standards}\\
Funding: This work was partially supported by the Istituto Nazionale di Alta Matematica - Gruppo Nazionale per il Calcolo Scientifico (INdAM-GNCS), Italy.
{Open Access funding enabled and organized by Italy Transformative Agreement.}
\\
Conflict of Interest:  The authors declare that they have no conflict of interest.
\\
Ethical approval: This article does not contain any studies with human participants or animals performed by any of the authors.
\\
\textbf{Acknowledgments} Francesco Calabr\`o  and Maria Roberta Belardo are members of GNCS-INdAM.  The authors would like to express many thanks to anonymous referees for helpful corrections and valuable comments on the original version of this paper.\B 
\\
\textbf{Data availability}
 The instructions for generating the datasets supporting the conclusions of this article are included in the article. Any other detail regarding the coding can be provided by the corresponding author upon reasonable request. \B
\\
\textbf{Open Access} This article is licensed under a Creative Commons Attribution 4.0 International License, which
permits use, sharing, adaptation, distribution and reproduction in any medium or format, as long as you give
appropriate credit to the original author(s) and the source, provide a link to the Creative Commons licence,
and indicate if changes were made. The images or other third party material in this article are included in the
article’s Creative Commons licence, unless indicated otherwise in a credit line to the material. If material is
not included in the article’s Creative Commons licence and your intended use is not permitted by statutory
regulation or exceeds the permitted use, you will need to obtain permission directly from the copyright holder.
To view a copy of this licence, visit http://creativecommons.org/licenses/by/4.0/ .



\bibliographystyle{plain}
\bibliography{biblio.bib}

\begin{thebibliography}{10}

\bibitem{10.1093/imanum/dry024}
B.~Adcock, R.~B. Platte, and A.~Shadrin.
\newblock Optimal sampling rates for approximating analytic functions from
  pointwise samples.
\newblock {\em IMA Journal of Numerical Analysis}, 39(3):1360--1390, 05 2018.

\bibitem{barron1993universal}
A.~R. Barron.
\newblock Universal approximation bounds for superpositions of a sigmoidal
  function.
\newblock {\em IEEE Transactions on Information theory}, 39(3):930--945, 1993.

\bibitem{battles2004extension}
Z.~Battles and L.~N. Trefethen.
\newblock An extension of matlab to continuous functions and operators.
\newblock {\em SIAM Journal on Scientific Computing}, 25(5):1743--1770, 2004.

\bibitem{bellman1957dimensionality}
R.~E. Bellman.
\newblock {\em Dynamic programming}.
\newblock Princeton University Press, 1957.

\bibitem{bishop2006pattern}
C.~M. Bishop.
\newblock {\em Pattern recognition and machine learning}.
\newblock springer, 2006.

\bibitem{boyd2009exponentially}
J.~P. Boyd and J.~R. Ong.
\newblock Exponentially-convergent strategies for defeating the runge
  phenomenon for the approximation of non-periodic functions, part i:
  single-interval schemes.
\newblock {\em Comput. Phys}, 5(2-4):484--497, 2009.

\bibitem{broomhead1988radial}
D.~Broomhead and D.~Lowe.
\newblock Radial basis functions, multi-variable functional interpolation and
  adaptive networks.
\newblock {\em Royal Signals and Radar Establishment Malvern (UK)}, 4148, 03
  1988.

\bibitem{calabro2009evaluation}
F.~Calabr{\`o} and A.~C. Esposito.
\newblock An evaluation of clenshaw--curtis quadrature rule for integration wrt
  singular measures.
\newblock {\em Journal of computational and applied mathematics},
  229(1):120--128, 2009.

\bibitem{calabro2021extreme}
F.~Calabr{\`o}, G.~Fabiani, and C.~Siettos.
\newblock Extreme learning machine collocation for the numerical solution of
  elliptic pdes with sharp gradients.
\newblock {\em Computer Methods in Applied Mechanics and Engineering},
  387:114188, 2021.

\bibitem{corless2020runge}
R.~M. Corless and L.~R. Sevyeri.
\newblock The {R}unge example for interpolation and {W}ilkinson's examples for
  rootfinding.
\newblock {\em SIAM Review}, 62(1):231--243, 2020.

\bibitem{cyr2020robust}
E.~C. Cyr, M.~A. Gulian, R.~G. Patel, M.~Perego, and N.~A. Trask.
\newblock Robust training and initialization of deep neural networks: An
  adaptive basis viewpoint.
\newblock In {\em Mathematical and Scientific Machine Learning}, pages
  512--536. PMLR, 2020.

\bibitem{ding2015extreme}
S.~Ding, H.~Zhao, Y.~Zhang, X.~Xu, and R.~Nie.
\newblock Extreme learning machine: algorithm, theory and applications.
\newblock {\em Artificial Intelligence Review}, 44:103--115, 2015.

\bibitem{dong2022computing}
S.~Dong and J.~Yang.
\newblock On computing the hyperparameter of extreme learning machines:
  Algorithm and application to computational pdes, and comparison with
  classical and high-order finite elements.
\newblock {\em Journal of Computational Physics}, 463:111290, 2022.

\bibitem{driscoll2014chebfun}
T.~A. Driscoll, N.~Hale, and L.~N. Trefethen.
\newblock {\em Chebfun guide}.
\newblock Pafnuty Publications, Oxford, 2014.

\bibitem{weinan2020towards}
W.~E, C.~Ma, S.~Wojtowytsch, and L.~Wu.
\newblock Towards a mathematical understanding of neural network-based machine
  learning: What we know and what we don’t.
\newblock {\em arXiv:2009.10713}, 2020.

\bibitem{ma2022barron}
W.~E, C.~Ma, and L.~Wu.
\newblock The {B}arron space and the flow-induced function spaces for neural
  network models.
\newblock {\em Constructive Approximation}, 55(1):369--406, 2022.

\bibitem{fornasier2022finite}
M.~Fornasier, T.~Klock, M.~Mondelli, and M.~Rauchensteiner.
\newblock Finite sample identification of wide shallow neural networks with
  biases.
\newblock {\em arXiv preprint arXiv:2211.04589}, 2022.

\bibitem{fornberg2011stable}
B.~Fornberg, E.~Larsson, and N.~Flyer.
\newblock Stable computations with gaussian radial basis functions.
\newblock {\em SIAM Journal on Scientific Computing}, 33(2):869--892, 2011.

\bibitem{han2018solving}
J.~Han, A.~Jentzen, and E.~Weinan.
\newblock Solving high-dimensional partial differential equations using deep
  learning.
\newblock {\em Proceedings of the National Academy of Sciences},
  115(34):8505--8510, 2018.

\bibitem{higham2019deep}
C.~F. Higham and D.~J. Higham.
\newblock Deep learning: An introduction for applied mathematicians.
\newblock {\em SIAM Review}, 61(4):860--891, 2019.

\bibitem{hornik1990universal}
K.~Hornik, M.~Stinchcombe, and H.~White.
\newblock Universal approximation of an unknown mapping and its derivatives
  using multilayer feedforward networks.
\newblock {\em Neural networks}, 3(5):551--560, 1990.

\bibitem{hryniowski2019deeplabnet}
A.~Hryniowski and A.~Wong.
\newblock Deeplabnet: End-to-end learning of deep radial basis networks with
  fully learnable basis functions.
\newblock {\em arXiv preprint arXiv:1911.09257}, 2019.

\bibitem{huang2015trends}
G.~Huang, G.-B. Huang, S.~Song, and K.~You.
\newblock Trends in extreme learning machines: A review.
\newblock {\em Neural Networks}, 61:32--48, 2015.

\bibitem{huang2006extreme}
G.-B. Huang, Q.-Y. Zhu, and C.-K. Siew.
\newblock Extreme learning machine: theory and applications.
\newblock {\em Neurocomputing}, 70(1-3):489--501, 2006.

\bibitem{jagtap2021deep}
A.~D. Jagtap, Y.~Shin, K.~Kawaguchi, and G.~E. Karniadakis.
\newblock Deep kronecker neural networks: A general framework for neural
  networks with adaptive activation functions.
\newblock {\em arXiv preprint arXiv:2105.09513}, 2021.

\bibitem{jin2017deep}
K.~H. Jin, M.~T. McCann, E.~Froustey, and M.~Unser.
\newblock Deep convolutional neural network for inverse problems in imaging.
\newblock {\em IEEE Transactions on Image Processing}, 26(9):4509--4522, 2017.

\bibitem{karniadakis2021physics}
G.~E. Karniadakis, I.~G. Kevrekidis, L.~Lu, P.~Perdikaris, S.~Wang, and
  L.~Yang.
\newblock Physics-informed machine learning.
\newblock {\em Nature Reviews Physics}, pages 1--19, 2021.

\bibitem{kratsios2021universal}
A.~Kratsios.
\newblock The universal approximation property: Characterizations, existence,
  and a canonical topology for deep-learning.
\newblock {\em Annals of Mathematics and Artificial Intelligence},
  89(5-6):435--469, 2021.

\bibitem{leshno1993multilayer}
M.~Leshno, V.~Y. Lin, A.~Pinkus, and S.~Schocken.
\newblock Multilayer feedforward networks with a nonpolynomial activation
  function can approximate any function.
\newblock {\em Neural networks}, 6(6):861--867, 1993.

\bibitem{lu2021learning}
L.~Lu, P.~Jin, G.~Pang, Z.~Zhang, and G.~E. Karniadakis.
\newblock Learning nonlinear operators via deeponet based on the universal
  approximation theorem of operators.
\newblock {\em Nature Machine Intelligence}, 3(3):218--229, 2021.

\bibitem{Mis2021}
S.~Mishra and R.~Molinaro.
\newblock Estimates on the generalization error of physics-informed neural
  networks for approximating a class of inverse problems for {PDEs}.
\newblock {\em {IMA} Journal of Numerical Analysis}, jun 2021.

\bibitem{neufeld2023universal}
A.~Neufeld and P.~Schmocker.
\newblock Universal approximation property of random neural networks.
\newblock {\em arXiv preprint arXiv:2312.08410}, 2023.

\bibitem{park1991universal}
J.~Park and I.~W. Sandberg.
\newblock Universal approximation using radial-basis-function networks.
\newblock {\em Neural computation}, 3(2):246--257, 1991.

\bibitem{pinkus}
A.~Pinkus.
\newblock Approximation theory of the mlp model.
\newblock {\em Acta Numerica 1999: Volume 8}, 8:143--195, 1999.

\bibitem{pinkus2015ridge}
A.~Pinkus.
\newblock {\em Ridge functions}, volume 205.
\newblock Cambridge University Press, 2015.

\bibitem{platte2011impossibility}
R.B. Platte, L.N Trefethen, and A.B.J. Kuijlaars.
\newblock Impossibility of fast stable approximation of analytic functions from
  equispaced samples.
\newblock {\em SIAM review}, 53(2):308--318, 2011.

\bibitem{qu2016two}
B.~Qu, B.F. Lang, J.~J. Liang, A.~K. Qin, and O.~D. Crisalle.
\newblock Two-hidden-layer extreme learning machine for regression and
  classification.
\newblock {\em Neurocomputing}, 175:826--834, 2016.

\bibitem{siegel2020approximation}
J.~W. Siegel and J.~Xu.
\newblock Approximation rates for neural networks with general activation
  functions.
\newblock {\em Neural Networks}, 128:313--321, 2020.

\bibitem{siegel2022high}
J.~W. Siegel and J.~Xu.
\newblock High-order approximation rates for shallow neural networks with
  cosine and reluk activation functions.
\newblock {\em Applied and Computational Harmonic Analysis}, 58:1--26, 2022.

\bibitem{trefethen2008gauss}
L.~N. Trefethen.
\newblock Is gauss quadrature better than clenshaw--curtis?
\newblock {\em SIAM review}, 50(1):67--87, 2008.

\bibitem{trefethen2019approximation}
L.~N. Trefethen.
\newblock {\em Approximation Theory and Approximation Practice, Extended
  Edition}.
\newblock SIAM, 2019.

\bibitem{vidal2017mathematics}
R.~Vidal, J.~Bruna, R.~Giryes, and S.~Soatto.
\newblock Mathematics of deep learning.
\newblock {\em arXiv preprint arXiv:1712.04741}, 2017.

\bibitem{wang2022review}
J.~Wang, S.~Lu, S.-H. Wang, and Y.-D. Zhang.
\newblock A review on extreme learning machine.
\newblock {\em Multimedia Tools and Applications}, 81(29):41611--41660, 2022.

\bibitem{wang2011study}
Y.~Wang, F.~Cao, and Y.~Yuan.
\newblock A study on effectiveness of extreme learning machine.
\newblock {\em Neurocomputing}, 74(16):2483--2490, 2011.

\bibitem{yuan2011optimization}
Y.~Yuan, Y.~Wang, and F.~Cao.
\newblock Optimization approximation solution for regression problem based on
  extreme learning machine.
\newblock {\em Neurocomputing}, 74(16):2475--2482, 2011.

\end{thebibliography}

\end{document}